\documentclass[final,3p,times]{elsarticle}
\usepackage{amssymb}
\usepackage{amsmath}
 \usepackage{amsthm}

\usepackage[T1]{fontenc}
\usepackage{tikz, mdframed, subcaption, amsmath, multirow, tikz, graphicx, lmodern, enumitem}
\usetikzlibrary{tikzmark,calc}

\newtheorem{theorem}{Theorem}

\newtheorem{corollary}{Corollary}
\newtheorem{lemma}{Lemma}
\newtheorem{example}{Example}
\newtheorem{reduction}{Reduction}
\newtheorem{remark}{Remark}
\newtheorem{observation}{Observation}

{\bfseries}{\itshape}
{\bfseries}{\itshape}
{\bfseries}{\itshape}
{\bfseries}{\itshape}
{\bfseries}{\itshape}

\newcommand{\x}{\gamma^{\rm{X}}}
\newcommand{\id}{\gamma^{\rm{ID}}}
\newcommand{\itd}{\gamma^{\rm{ITD}}}
\newcommand{\fd}{\gamma^{\rm{FD}}}
\newcommand{\ftd}{\gamma^{\rm{FTD}}}
\newcommand{\ld}{\gamma^{\rm{LD}}}
\newcommand{\ltd}{\gamma^{\rm{LTD}}}
\newcommand{\od}{\gamma^{\rm{OD}}}
\newcommand{\otd}{\gamma^{\rm{OTD}}}

\newcommand{\X}{{\textsc{X}}}
\newcommand{\ID}{{\textsc{ID}}}
\newcommand{\ITD}{{\textsc{ITD}}}
\newcommand{\FD}{{\textsc{FD}}}
\newcommand{\FTD}{{\textsc{FTD}}}
\newcommand{\LD}{{\textsc{LD}}}
\newcommand{\LTD}{{\textsc{LTD}}}
\newcommand{\OD}{{\textsc{OD}}}
\newcommand{\OTD}{{\textsc{OTD}}}

\newcommand{\fdadmis}{{\textsc{FD}}-admissible}
\newcommand{\ftdadmis}{{\textsc{FTD}}-admissible}

\newcommand{\codes}{{\textsc{Codes}}}
\newcommand{\hyp}{\mathcal{H}}

\newcommand{\for}{{\textsc{For}}}
\newcommand{\calC}{\mathcal{C}}
\newcommand{\NP}{$\mathsf{NP}$}
\newcommand{\yes}{$\mathsf{YES}$}

\newcommand{\mincode}[1]{\textsc{Min #1-Code}}

\newcommand{\ddelta}{\bigtriangleup}

\newcommand{\F}{\mathcal{F}}

\newcommand{\calO}{\mathcal{O}}

\newtheorem{claim}{$\blacksquare$ Claim}[lemma]

\newenvironment{claimproof}{\noindent\emph{Proof of claim.}}{\hfill$\blacksquare$\newline\medskip}

\newcommand{\defproblem}[3]{
\vspace{3mm}
\noindent \fbox{
\begin{minipage}{0.96\textwidth}
\begin{tabular*}{\textwidth}{@{\extracolsep{\fill}}lr}
#1 \\ \end{tabular*}
{\bf{Input:}} #2 \\
{\bf{Question:}} #3
\end{minipage}
}

\vspace{3mm}
}

\newcommand{\w}[1]{{\color{white} #1}}
\newcommand{\bw}[1]{{\color{white}{\textbf{#1}}}}
\newcommand{\bt}[1]{{\textbf{#1}}}

\begin{document}

\begin{frontmatter}



\title{On full-separating sets and related codes in graphs\footnote{We dedicate this work to Ralf Klasing to thank him for inspiring and accompanying our work during several meetings and to express our deep appreciation for his scientific work.}}


\author[label1,label2]{Dipayan Chakraborty}
\ead{dipayancha@gmail.com}%% Author name
\author[label1]{Annegret K. Wagler} 
\ead{annegret.wagler@uca.fr}%% Author name

\affiliation[label1]{organization={Universit\'{e} Clermont-Auvergne, CNRS, Mines de Saint-\'{E}tienne, Clermont-Auvergne-INP, LIMOS},%Department and Organization
            %addressline={} 
            city={Clermont-Ferrand},
            postcode={63000}, 
            %state={}
            country={France}}
            
\affiliation[label2]{organization={Department of Mathematics and Applied Mathematics, University of Johannesburg},%Department and Organization
            %addressline={}, 
            city={Auckland Park},
            postcode={2006}, 
            %state={}
            country={South Africa}}

\begin{abstract}
A domination-based identification problem on a graph $G$ is one where the objective is to choose a subset $C$ of the vertex set of $G$ such that $C$ has both, a \emph{domination property}, that is, $C$ is either a dominating or a total-dominating set of $G$, and a \emph{separation property}, that is, any two distinct vertices of $G$ must have distinct closed or open neighborhoods in $C$. Such a set $C$ is often referred to as a \emph{code} in the literature of identification problems. In this article, we introduce a new separation property, called \emph{full-separation}, as it combines aspects of the two well-studied properties of closed- and open-separation. We study it in combination with both domination and total-domination and call the resulting codes \emph{full-separating dominating codes} (or \emph{FD-codes} for short) and \emph{full-separating total-dominating codes} (or \emph{FTD-codes} for short), respectively. Incidentally, FTD-codes have also been introduced in the literature of identification problems under the name of \emph{strongly identifying codes} (or \emph{SID-codes} for short) and under a differently formulated definition. In this paper, we address the conditions for the existence of FD- and FTD-codes, bounds for their size, their relation to codes of the other types and present some extremal cases for these bounds and relations. We further show that the problems of determining an FD- or an FTD-code of minimum cardinality in a graph are NP-hard. We also show that the cardinalities of minimum FD- and FTD-codes of any graph differ by at most one, but that it is NP-hard to decide whether or not they are equal for a given graph in general. 
\end{abstract}


\end{frontmatter}

\emph{Ralf Klasing introduced me, Annegret Wagler, to the area of identification problems in graphs on the occasion of a workshop held in Bordeaux in 2011. This workshop inspired me to start working in this area. Being familiar with polyhedral combinatorics, I wondered how this approach can be applied to identification problems. In collaboration with some experts on covering polyhedra from Argentina, we studied identification problems from this point of view, based on their reformulation in terms of covering problems in suitably constructed hypergraphs. Analyzing these hypergraph representations, we, Dipayan Chakraborty and myself, discovered full-separation as new separation property and show here how this opens the perspective for a systematic study of identification problems.}

\section{Introduction}

\begin{figure}[!t]
\begin{center}
\includegraphics[scale=0.9]{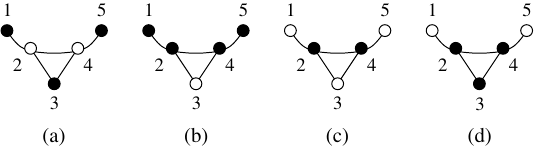}
\caption{Minimum X-codes in the \emph{bull} graph (the black vertices belong to the code), where (a) is an ID-code, (b) an ITD-code, (c) both an LD- and LTD-code, (d) both an OD- and OTD-code.}
\label{fig_exp_X-codes}
\end{center}
\end{figure}

In the domain of identification problems on graphs, two distinct vertices of a graph $G$ are said to be \emph{separated} if they have distinct neighborhoods in a suitably chosen dominating or total-dominating set of $G$. Depending on the types of dominating sets and separation properties used, various problems arise under various names in the literature. More precisely, consider a graph $G=(V,E)$ and denote by $N(v) = \{u \in V : uv \in E\}$ (respectively, $N[v] = N(v) \cup \{v\}$) the open (respectively, closed) neighborhood of a vertex $v \in V$. A subset $C \subseteq V$ is \emph{dominating} (respectively, \emph{total-dominating}) if the set $N[v]\cap C$ (respectively, $N(v)\cap C$) is non-empty for each $v \in V$. On the other hand, a subset $C \subseteq V$ is \emph{closed-separating} (respectively, \emph{open-separating}) if the set $N[v]\cap C$ (respectively, $N(v)\cap C$) is unique for each $v \in V$. Moreover, $C$ is \emph{locating} if the set $N(v)\cap C$ is unique for each $v \in V\setminus C$. So far, the following combinations of separation and domination properties have been studied in the literature. \\[-5mm]
\begin{itemize}[leftmargin=12pt, itemsep=0pt]
%%%%%%%%%
\item closed-separation with domination and total-domination leading to \emph{identifying codes} (\emph{ID-codes} for short) and \emph{identifying total-dominating codes}\footnote{ITD-codes had been introduced in~\cite{HHH_2006} under the name \emph{differentiating total-do\-mi\-na\-ting codes}. Due to consistency in notation, we call them ITD-codes in this article.} (\emph{ITD-codes} for short), see \cite{KCL_1998} and \cite{HHH_2006}, respectively;
%%%%%%%%%%%
\item location with domination and total-domination leading to \emph{locating domina\-ting codes} (\emph{LD-codes} for short) and \emph{locating total-dominating codes} (\emph{LTD-codes} for short), see \cite{S_1988} and \cite{HHH_2006}, respectively;
%%%%%%%%%%%%
\item open-separation with domination and total-domination leading to \emph{open-sepa\-ra\-ting dominating codes} (\emph{OD-codes} for short) and \emph{open-separating total-do\-minating codes}
\footnote{OTD-codes were introduced independently in~\cite{HLR_2002} and in~\cite{SS_2010} under the names of \emph{identifying code with nontransmitting vertices} (or \emph{IDNT-code}) and \emph{open neighborhood locating-dominating sets} (or \emph{OLD-sets}), respectively. However, due to consistency in naming that specifies the separation and the domination property, we prefer to call them open-separating total-dominating codes in this article.
} (\emph{OTD-codes} for short), see \cite{CW_ISCO2024} and  \cite{HLR_2002,SS_2010}, respectively. \\[-5mm]
\end{itemize}
%%%%%%%%%%%%%%
Figure \ref{fig_exp_X-codes} illustrates examples of such codes in a small graph. 
%%%%%%%%%%%%%%
%%%%%%%%%%%%%%
Given a graph $G=(V,E)$, for X $\in$ \{ID, ITD, LD, LTD, OD, OTD\}, the X-problem on $G$ is the problem of finding an X-code of minimum cardinality $\x(G)$ in $G$. The parameter $\x(G)$ is called the \emph{\X-number} of the graph $G$. 

Problems of this type are intensively studied in the literature (see, for example, the internet bibliography containing over 500 articles around these topics maintained by Jean and Lobstein~\cite{JL_lib}). Such problems have manifold applications, for example, in locating intruders in facilities by placing detectors (with sensors) in some of the rooms such that they can detect intruders up to their neighboring rooms~\cite{UTS04}. Similar models apply to detecting faults in a network of multiprocessors where some of the multiprocessors are designated as detectors and can detect faulty multiprocessors which are their neighbors in the network. All these practical scenarios can be modeled as a graph-theoretic problem where a graph represents a network/ facility, its vertices represent multiprocessors/ rooms and edges represent adjacencies between them. In other words, the problem translates to choosing the code vertices (preferably a ``small'' subset of the vertex set) that represent detectors in the network/ facility. 

Each of the above separation properties of location, closed-separation and open-separation, models specific requirements of fault detection. For example, imagine the detectors used are ``3-state'' detectors, that is, each detector sends the signal~0 if no fault occurs in its closed neighborhood, it sends the signal~$1$ if a fault occurs at a neighboring vertex of the detector and sends a signal~$2$ if the fault occurs at the same vertex as the detector. In other words, such detectors have the capacity to distinguish between ``no fault in the neighborhood'', ``fault at a neighbor'' and ``fault at itself''. In such a case, distributing these detectors in the network in the form of a locating set~\cite{S_1988} effectively identifies faults in the network by assigning a unique signal string to each vertex $v$. Assigning the signal string may be done the following way: Let $v_1,v_2, \ldots , v_n$ be the vertices of a graph $G$. Then assign to each vertex $v$ the signal string $f(v) = x_1x_2 \ldots x_n$, where 
\begin{eqnarray*}
x_i = \begin{cases}
        0, &\text{ if $v_i$ has no detector on it};\\
        \text{signal value emitted by $v_i$ when there is a fault at $v$}, &\text{ if $v_i$ has a detector}.
        \end{cases}
\end{eqnarray*}
For example, let $v_1,v_2,v_3,v_4,v_5$ be the vertices $1,2,3,4,5$, respectively, of the \emph{bull} graph shown in Figure~\ref{fig_exp_X-codes}. Then, as shown by means of the black vertices in Figure~\ref{fig_exp_X-codes}(c), distributing the 3-state detectors on the two black vertices (which is a locating set of the graph), the vertices $v_1,v_2,v_3,v_4,v_5$ are assigned the (unique) signal strings $01000$, $02010$, $01010$, $01020$, $00010$, respectively. Notice that the presence of $x_i=2$ in a signal string immediately implies the presence of the fault at the vertex $v_i$. For example, $x_4=2$ in $f(v_4) = 01020$ of Figure~\ref{fig_exp_X-codes}(c) implies that the fault definitively occurs at the vertex $v_4$. Hence, it is only the set of signal strings with all $x_i \in \{0,1\}$ that have to be pairwise distinct.

Let us now assume that the detectors used to monitor a network are ``2-state'' detectors (of Type~1), which send the signal~1 if a fault occurs anywhere in its closed neighborhood and~0 otherwise. This is the same as replacing all the 2's in the signal strings produced by the 3-state detectors by~1's. In other words, these detectors lack the ability to distinguish between ``fault at a neighbor'' and ``fault at itself''. With such detectors, the network monitoring modeled by locating sets fails. For example, if the detectors at vertices $v_2$ and $v_4$ of Figure~\ref{fig_exp_X-codes}(c) are 2-state of Type~1, then their respective signal strings $f(v_2)$ and $f(v_4)$ both equal $01010$ (that is, replacing the 2's in the previous signal strings by 1's). Hence, with a locating set, a fault at any of the vertices $v_2$ and $v_4$ is not detected. Therefore, to effectively monitor networks with 2-state Type~1 detectors, one needs the property of closed-separation~\cite{KCL_1998}. For example, distributing the 2-state Type~1 detectors on the black vertices in Figure~\ref{fig_exp_X-codes}(a) (which represents a closed-separating set), the vertices $v_1,v_2,v_3,v_4,v_5$ have their (unique) %codes
signal strings $10000$, $10100$, $00100$, $00101$, $00001$, respectively.

The third variation mentioned above, namely the concept of open-separating sets~\cite{HLR_2002,SS_2010}, can be seen as those modeling network-monitoring by ``2-state'' detectors (of Type 2), whereby such detectors send the signal 1 if a fault occurs at a neighboring vertex and 0 otherwise (which includes its own vertex). This is the same as replacing all the 2's in the signal strings produced by the 3-state detectors by~0's. In this case too, distributing such detectors in the form of a locating set does not ensure effective monitoring of a network. For example, with such 2-state Type~2 detectors, the signal strings $f(v_2)$ and $f(v_5)$ in Figure~\ref{fig_exp_X-codes}(c) both equal $00010$ (that is, by replacing the 2's in signal strings of the 3-state detectors by 0's) and hence, a fault at any of the vertices $v_2$ and $v_5$ is not detected. However, distributing these detectors in the form of an open-separating set, for example, on the black vertices in Figure~\ref{fig_exp_X-codes}(d) (which represents an open-separating set), yields the (unique) signal strings $01000$, $00110$, $01010$, $01100$, $00010$ for the vertices $v_1,v_2,v_3,v_4,v_5$, respectively.

The purpose of this article is to introduce a separation property, called \emph{full-separation}, whose practical motivation lies in modeling a monitoring system in which the 3-state detectors may themselves suffer from errors causing them to be transformed into 2-state detectors, either of Type~1 or Type~2 that is, the signals 2's from some of the 3-state detectors may be replaced by either~1's or~0's). Such a model requires both aspects of closed-separation and open-separation. As we shall also see later in analyzing the patterns in Table~\ref{tab_hypergraphs}, there is also a theoretical motivation to introduce the property of full-separation. We formally introduce full-separation later in Section~\ref{sec2}.

For $\X \in \{\LD, \LTD, \ID, \ITD, \OD, \OTD \}$, the X-problems have been studied from a unifying point of view, namely as reformulations in terms of covering problems in suitably constructed hypergraphs, for example, in \cite{ABLW_2018,ABLW_2022,ABW_2016,CW_ISCO2024}. 
Given a graph $G=(V,E)$ and an X-problem, we look for a hypergraph $\hyp_\X(G) = (V,\F_\X)$ so that $C \subseteq V$ is an X-code of $G$ if and only if $C$ is a \emph{cover} of $\hyp_\X(G)$, that is, a subset $C \subseteq V$ satisfying $C \cap F \neq \emptyset$ for all $F \in \F_\X$. Then the \emph{covering number} $\tau(\hyp_\X(G))$, defined as the minimum cardinality of a cover of $\hyp_\X(G)$, equals by construction the X-number $\x(G)$. The hypergraph $\hyp_\X(G)$ is called the \emph{\X-hypergraph} of $G$. Refer to Table~\ref{tab_exp_Hyps} for examples of such \X-hypergraphs for the bull graph shown in Figure~\ref{fig_exp_X-codes}.

It is a simple observation that for an X-problem involving domination (respectively, total-domination), $\F_\X$ needs to contain the closed (respectively, open) neighborhoods of all vertices of $G$. In order to encode the separation properties,  that is, the fact that the intersection of an X-code with the neighborhood of each vertex is \emph{unique}, it was suggested in \cite{ABLW_2018,ABLW_2022} to use the symmetric differences of the neighborhoods. Here, given two sets $A$ and $B$, their \emph{symmetric difference} is defined by $A \bigtriangleup B = (A \setminus B) \cup (B \setminus A)$. In fact, it has been shown in \cite{ABLW_2018,ABLW_2022} that a code $C$ of a graph $G$ is  \\[-5mm]
%%%%%%%%%%%%%%%%%%
\begin{itemize}[leftmargin=12pt, itemsep=0pt]
\item closed-separating if and only if $(N[u] \bigtriangleup N[v]) \cap C \neq \emptyset$ for all distinct $u,v \in V$,
\item open-separating if and only if $(N(u) \bigtriangleup N(v)) \cap C \neq \emptyset$ for all distinct $u,v \in V$, 
\item locating if and only if $(N(u) \bigtriangleup N(v)) \cap C \neq \emptyset$ for all $uv \in E$ and $(N[u] \bigtriangleup N[v]) \cap C \neq \emptyset$ for all $uv \notin E$, where $u,v \in V$ are distinct. \\[-5mm]
\end{itemize}
%%%%%%%%%%%%%%%%%%
Accordingly, for each of the already studied X-problems, the hyperedge set $\F_\X$ of $\hyp_\X(G)$ consists of exactly 3 different subsets as shown in Table~\ref{tab_hypergraphs}, where \\[-5mm]
%%%%%%%%%%%%%%%%%
\begin{itemize}[leftmargin=12pt, itemsep=0pt]
\item $N[G]$ (respectively, $N(G)$) is the set of all closed (respectively, open) neighborhoods of vertices in $G$, 
\item $\bigtriangleup_a[G]$ (respectively, $\bigtriangleup_a(G)$) denotes the set of symmetric differences of closed (respectively, open) neighborhoods of all pairs of \emph{adjacent} vertices in $G$, 
\item $\bigtriangleup_n[G]$ (respectively, $\bigtriangleup_n(G)$) denotes the set of symmetric differences of closed (respectively, open) neighborhoods of all pairs of \emph{non-adjacent} vertices in $G$. \\[-5mm]
\end{itemize}
%%%%%%%%%%%%%%%%%%
\begin{table}[t]
\begin{center}
\begin{tabular}{ c || c | c || c | c || c | c }
 X & ID   & ITD  & LD   & LTD  & OD & OTD \\ \hline 
   & $N[G]$ & $N(G)$ & $N[G]$ & $N(G)$ & $N[G]$ & $N(G)$ \\
$\F_\X$ & $\bigtriangleup_a[G]$~ & ~$\bigtriangleup_a[G]$~ & ~$\bigtriangleup_a(G)$~ & ~$\bigtriangleup_a(G)$~ & ~$\bigtriangleup_a(G)$~ & ~$\bigtriangleup_a(G)$ \\
 &  ~$\bigtriangleup_n[G]$~ & ~$\bigtriangleup_n[G]$~ & ~$\bigtriangleup_n[G]$~ & ~$\bigtriangleup_n[G]$~ & ~$\bigtriangleup_n(G)$~ & ~$\bigtriangleup_n(G)$ \\
\end{tabular}
\end{center}
\caption{
The X-hypergraphs $\hyp_\X(G) = (V,\F_\X)$ for the already considered X-problems, listing column-wise the 3 needed subsets of hyperedges}
\label{tab_hypergraphs}
\end{table}
%%%%%%%%%%%%%%%%%%%
For illustration, Table \ref{tab_exp_Hyps} shows the resulting \X-hypergraphs for the bull (Figure~\ref{fig_exp_X-codes}) and four different X-codes, namely $\X \in \{\ID, \ITD, \LTD, \OTD\}$.
%%%%%%%%%%%%%%%%%%
\begin{table}[ht]
\begin{center}
\begin{tabular}{ c | c | c | c | c  }
    & ID             & ITD  & LTD  &  OTD \\ \hline\hline
$v \in V(G)$ & $N[v]$ & $N(v)$ & $N(v)$ & $N(v)$ \\ \hline
  1 & \{ \bt{1}, \ 2~ \ \bw{3,} \ \w{4,} \ \bw{5,} \} & \{ \bw{1,} \ \bt{2}~ \ \w{3,} \ \bw{4,} \ \bw{5,} \} & \{ \w{1,} \ \bt{2}~ \ \w{3,} \ \bw{4,} \ \w{5,} \} & \{ \w{1,} \ \bt{2}~ \ \bw{3,} \ \bw{4,} \ \w{5,} \} \\
  2 & \{ \bt{1}, \ 2, \ \bt{3}, \ 4~ \ \bw{5,} \} & \{ \bt{1}, \ \bw{2,} \ 3, \ \bt{4}, \ \bw{5,} \} & \{ 1, \ \bw{2,} \ 3, \ \bt{4}, \ \w{5,} \}                & \{ 1,     \ \bw{2,} \ \bt{3}, \ \bt{4}~ \ \w{5,} \} \\
  3 & \{ \bw{1,} \ 2, \ \bt{3}, \ 4~ \ \bw{5,} \} & \{ \bw{1,} \ \bt{2}, \ \w{3,} \ \bt{4}~ \ \bw{5,} \} & \{ \w{1,} \ \bt{2}, \ \w{3,} \ \bt{4}~ \ \w{5,} \}    & \{ \w{1,} \ \bt{2}, \ \bw{3,} \ \bt{4}~ \ \w{5,} \} \\
  4 & \{ \bw{1,} \ 2, \ \bt{3}, \ 4, \ \bt{5}~ \} & \{ \bw{1,} \ \bt{2}, \ 3, \ \bw{4,} \ \bt{5}~ \} & \{ \w{1,} \ \bt{2}, \ 3, \ \bw{4,} \ 5~ \}                & \{ \w{1,} \ \bt{2}, \ \bt{3}, \ \bw{4,} \ 5~ \} \\
  5 & \{ \bw{1,} \ \w{2,} \ \bw{3,} \ 4, \ \bt{5}~ \} & \{ \bw{1,} \ \bw{2,} \ \w{3,} \ \bt{4}~ \ \bw{5}~ \} & \{ \w{1,} \ \bw{2,} \ \w{3,} \ \bt{4}~ \ \w{5}~ \} & \{ \w{1,} \ \bw{2,} \ \bw{3,} \ \bt{4}~ \ \w{5}~ \} \\ \hline\hline
$uv \in E(G)$ & $N[u] \bigtriangleup N[v]$ & $N[u] \bigtriangleup N[v]$ & $N(u) \bigtriangleup N(v)$ & $N(u) \bigtriangleup N(v)$ \\ \hline
1,2 & \{ \bw{1,} \ \w{2,} \ \bt{3}, \ 4~ \ \bw{5,} \} & \{ \bw{1,} \ \w{2,} \ 3, \ \bt{4}~ \ \bw{5,} \} & \{ 1, \ \bt{2}, \ 3, \ \bt{4}~ \ \w{5,} \} & \{ 1, \ \bt{2}, \ \bt{3}, \ \bt{4}~ \ \w{5,} \} \\
2,3 & \{ \bt{1}~ \ \w{2,} \ \bw{3,} \ \w{4,} \ \bw{5,} \} & \{ \bt{1}~ \ \w{2,} \ \bw{3,} \ \w{4,} \ \bw{5,} \} & \{ 1, \ \bt{2}, \ 3~ \ \bw{4,} \ \w{5,} \} & \{ 1, \ \bt{2}, \ \bt{3}~ \ \bw{4,} \ \w{5,} \} \\
2,4 & \{ \bt{1}, \ \w{2,} \ \bw{3,} \ \w{4,} \ \bt{5}~ \} & \{ \bt{1}, \ \w{2,} \ \bw{3,} \ \w{4,} \ \bt{5}~ \} & \{ 1, \ \bt{2}, \ \w{3,} \ \bt{4}, \ 5~ \} & \{ 1, \ \bt{2}, \ \w{3,} \ \bt{4}, \ 5~ \} \\
3,4 & \{ \bw{1,} \ \w{2,} \ \bw{3,} \ \w{4,} \ \bt{5}~ \} & \{ \bw{1,} \ \w{2,} \ \bw{3,} \ \w{4,} \ \bt{5}~ \} & \{ \w{1,} \ \bw{2,} \ 3, \ \bt{4}, \ 5~ \} & \{ \w{1,} \ \bw{2,} \ \bt{3}, \ \bt{4}, \ 5~ \}\\
4,5 & \{ \bw{1,} \ 2, \ \bt{3}~ \ \w{4,} \ \bw{5,} \} & \{ \bw{1,} \ \bt{2}, \ 3~ \ \bw{4,} \ \bw{5,} \} & \{ \w{1,} \ \bt{2}, \ 3, \ \bt{4}, \ 5~ \} & \{ \w{1,} \ \bt{2}, \ \bt{3}, \ \bt{4}, \ 5~ \} \\ \hline\hline
$uv \notin V(G)$ & $N[u] \bigtriangleup N[v]$ & $N[u] \bigtriangleup N[v]$ & $N[u] \bigtriangleup N[v]$ & $N(u) \bigtriangleup N(v)$ \\ \hline
1,3 & \{ \bt{1}, \ \w{2,} \ \bt{3}, \ 4~ \ \bw{5,} \} & \{ \bt{1,} \ \bw{2,} \ 3, \ \bt{4}~ \ \bw{5,} \} & \{ 1, \ \bw{2,} \ 3, \ \bt{4}~ \ \w{5,} \} & \{ \w{1,} \ \bw{2,} \ \bw{3,} \ \bt{4}~ \ \w{5,} \} \\
1,4 & \{ \bt{1}, \ \w{2,} \ \bt{3}, \ 4, \ \bt{5}~ \} & \{ \bt{1}, \ \bw{2,} \ 3, \ \bt{4}, \ \bt{5}~ \} & \{ 1, \ \bw{2,} \ 3, \ \bt{4}, \ 5~ \} & \{ \w{1,} \ \bw{2,} \ \bt{3}, \ \bw{4,} \ 5~ \} \\
1,5 & \{ \bt{1}, \ 2, \ \bw{3,} \ 4, \ \bt{5}~ \} & \{ \bt{1}, \ \bt{2}, \ \w{3,} \ \bt{4}, \ \bt{5}~ \} & \{ 1, \ \bt{2}, \ \w{3,} \ \bt{4}, \ 5~ \} & \{ \w{1,} \ \bt{2}, \ \bw{3,} \ \bt{4}~ \ \w{5,} \} \\
2,5 & \{ \bt{1}, \ 2, \ \bt{3}, \ \w{4,} \ \bt{5}~ \} & \{ \bt{1}, \ \bt{2}, \ 3, \ \bw{4,} \ \bt{5}~ \} & \{ 1, \ \bt{2}, \ 3, \ \bw{4,} \ 5~ \} & \{ 1, \ \bw{2,} \ \bt{3}~ \ \bw{4,} \ \w{5,} \} \\
3,5 & \{ \bw{1,} \ 2, \ \bt{3}, \ \w{4,} \ \bt{5}~ \} & \{ \bw{1,} \ \bt{2}, \ 3, \ \bw{4,} \ \bt{5}~\} & \{ \w{1,} \ \bt{2}, \ 3, \ \bw{4,} \ 5~ \} & \{ \w{1,}  \ \bt{2}, \ \bw{3,} \ \bw{4,} \ \w{5}~ \} \\
\end{tabular}
\end{center}
\caption{
Illustration of the 3 needed subsets of hyperedges (neighborhoods of all vertices, symmetric differences of neighborhoods of pairs of adjacent and pairs of non-adjacent vertices) of the X-hypergraphs $\hyp_\X(G)$ where $G$ is the bull (see Figure \ref{fig_exp_X-codes} for the graph and the vertex numbers); and (a) $\X=\ID$, (b) $\X=\ITD$, (c) $\X=\LTD$, (d) $\X=\OTD$. The vertex numbers in bold exhibit the covers of $\hyp_\X(G)$ corresponding to the X-codes of $G$ shown by means of the black vertices in Figure \ref{fig_exp_X-codes}.}
\label{tab_exp_Hyps}
\end{table}
%%%%%%%%%%%%%%%%%%%

Analyzing the summary of hypergraph representations of the six already established X-problems in Table~\ref{tab_hypergraphs}, we observe that one separation property has not yet been considered in the literature, namely the one involving $\bigtriangleup_a[G]$ and $\bigtriangleup_n(G)$. The aim of this paper is to introduce the corresponding separation property, called \emph{full-separation}, whose hypergraph representations involve exactly the hyperedge subsets $\bigtriangleup_a[G]$ and $\bigtriangleup_n(G)$ missing from Table~\ref{tab_hypergraphs}. As we shall also see later in Section~\ref{sec2} where we introduce these concepts more formally, it turns out that a set is \emph{full-separating} if and only if it is both closed-separating and open-separating. Moreover, we study full-separation in combination with both domination and total-domination leading to \emph{full-separating dominating codes} (or \emph{FD-codes} for short) and \emph{full-separating total-dominating codes} (or \emph{FTD-codes} for short), respectively. As with the other existing codes in the literature, given a graph $G$, the parameters $\fd(G)$ and $\ftd(G)$, called the \emph{\FD-number} and the \emph{\FTD-number}, respectively, are the smallest cardinalities of \FD- and \FTD-codes, respectively, of $G$.

It is to be noted here that the notion of FTD-codes does turn out to be equivalent to the already existing notion of \emph{strongly identifying codes} (or \emph{SID-codes}) defined in~\cite{H_2010} (see Theorem~\ref{thm_PL-char} for a proof of this equivalence) and to the restriction $(\ell,t) = (1,1)$ on \emph{strongly $(\ell,\le t)$-identifying codes} defined on Hamming spaces in~\cite{HLR_2002}. Nevertheless, to the best of our knowledge, all further work on SID-codes since its introduction have been on either Hamming spaces or on grids (see for example,~\cite{H_2010, HLR_2001, HLR_2002, L_2002, R_2003}). In this paper, however, we initiate the study of FTD-codes (along with FD-codes) on graphs in general and under the umbrella of the full-separating property. Moreover, we reiterate that, apart from the above-mentioned practical motivations for the full-separating property, we also introduce the latter from the theoretical perspective of studying the hyperedge subsets $\bigtriangleup_a[G]$ and $\bigtriangleup_n(G)$ missing from Table~\ref{tab_hypergraphs}. Despite the already existing name of ``SID-code'', in this paper, we continue to refer to the said code as ``FTD'' since for most of this current work, we treat both FD- and FTD-codes together from the point of view of full-separation.

Note that not all graphs admit codes of all studied types. We accordingly address the conditions for the existence of FD- and FTD-codes, their relation to codes of the other types and bounds on $\fd(G)$ and $\ftd(G)$ in Section \ref{sec2}. In particular, we show that $\fd(G)$ and $\ftd(G)$ differ by at most one. Moreover, the problems of determining an X-code of minimum cardinality $\x(G)$ in a graph $G$ have been shown to be NP-hard for all the previously studied X-problems~\cite{CW_ISCO2024,CHL_2003,CSS_1987,SS_2010}. We show the same for FD- and FTD-codes in Section \ref{sec3}. Furthermore, we show that, despite the minimum cardinalities of the \FD- and the \FTD-codes differing by at most~$1$, it is NP-hard to decide whether or not the two said minimum cardinalities actually differ on a given graph in general. On a positive note, in Section~\ref{sec4}, we present closed formulas for the  FD- and FTD-numbers on some graphs classes like paths, cycles and some selected subclasses of bipartite and split graphs. Some of these examples in Section~\ref{sec4} also act as extremal cases for the \FD- and \FTD-numbers of graphs. Evidently therefore, also deciding if $\fd(G) \ne \ftd(G)$ for $G$ being in any of these latter graph classes can also be done in polynomial-time. Another interesting property of FD-codes that comes through in studying the said graph classes is that the FD-number of a graph can be strictly greater than the sum of the FD-numbers of its components. This phenomenon is not true for all X-codes discussed so far, except when X = OD~\cite{CW_ISCO2024}.

Some of the results presented in this article have also appeared in~\cite{CW_CALDAM2025} without proofs in the form of the conference version of this paper. For all basic definitions and notations on graph theory, we follow the book~\cite{West_2001} by West.

%%%%%%%%%%%%%%%%%%%%%%%%%%%%%%%%%%%%%%%%%%%

\section{Full-separation and related codes}
\label{sec2}

In this section, we formally introduce the concept of full-separation, the related full-separating codes and address fundamental questions concerning the existence of these codes and bounds for their cardinalities.

\subsection{Full-separation and related codes} 

\begin{figure}[h]
\begin{center}
\includegraphics[scale=0.8]{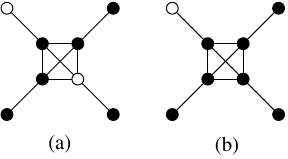}
\caption{Minimum X-codes in a graph (the black vertices belong to the code), where (a) is an FD-code, (b) an FTD-code.}
\label{fig_exp_PL-codes}
\end{center}
\end{figure}

Let $G=(V,E)$ be a graph. A subset $C \subseteq V$ is a \emph{full-separating set} of $G$ if, for each pair of distinct vertices $u,v \in V$, we have \\[-2mm]
\[
(N(v)\cap C) \setminus \{u\} \neq (N(u)\cap C) \setminus \{v\}.
\]
Moreover, a full-separating set $C \subseteq V$ is a \emph{full-separating dominating code} (\emph{FD-code} for short) (respectively, a \emph{full-separating total-dominating code} (\emph{FTD-code} for short)) if $C$ is also a dominating (respectively, total-dominating) set of $G$.  
For example, the ITD-code from Figure \ref{fig_exp_X-codes} is also both a minimum FD- and FTD-code of the graph. 
Moreover, Figure \ref{fig_exp_PL-codes} illustrates further examples of such \FD- and \FTD-codes in a small graph.
%%%%%%%%%%%%%%%
%%%%%%%%%%%%%%%%

\begin{observation}\label{apobs_nbhd_idty}
Let $G = (V,E)$ be a graph and $u,v$ be any two distinct vertices of $G$. Then, we have
\[
(N(u) \bigtriangleup N(v)) \setminus \{u,v\} = \left
\{ \begin{array}{lll}
N[u] \bigtriangleup N[v] &\subset N(u) \bigtriangleup N(v), &\text{ if $uv \in E$};\\
N(u) \bigtriangleup N(v) &\subset N[u] \bigtriangleup N[v], &\text{ if $uv \notin E$.}
\end{array}
\right.
\]
In particular, therefore, the set $(N(u) \bigtriangleup N(v)) \setminus \{u,v\}$ is a subset of both $N(u) \bigtriangleup N(v)$ and $N[u] \bigtriangleup N[v]$ for all distinct $u,v \in V$.
\end{observation}

In the following, we show some equivalent characterizations of full-separating sets. 
%%%%%%%%%%%%%%%%
\begin{theorem}\label{thm_PL-char}
Let $G=(V,E)$ be a graph. For a subset $C \subseteq V$, the following assertions are equivalent. \\[-5mm]
\begin{enumerate}[leftmargin=18pt, itemsep=0pt]
\item[(a)] \label{thm_PL-char_PL}
  $C$ is full-separating.
\item[(b)] \label{thm_PL-char_altdef}
  $C$ has a non-empty intersection with $\big( N(u) \bigtriangleup N(v) \big) \setminus \{u,v\}$ for all distinct $u,v \in V$.
\item[(c)] \label{thm_PL-char_adjacency}
  $C$ has a non-empty intersection with
\begin{itemize}[leftmargin=12pt, itemsep=0pt]
\item $N[u] \bigtriangleup N[v]$ for all pairs of adjacent vertices $u,v \in V$, and 
\item $N(u) \bigtriangleup N(v)$ for all pairs of non-adjacent vertices $u,v \in V$.
\end{itemize}
\item[(d)] \label{thm_PL-char_open-closed}
  $N[u] \cap C \ne N[v] \cap C$ and $N(u) \cap C \ne N(v) \cap C$ for each pair of distinct vertices $u,v \in V$.
\item[(e)] \label{thm_PL-char_open-closed_statement}
  $C$ is both a closed-separating and an open-separating set of $G$.
\item[(f)] \label{thm_PL-SID} The pairs $\{N(u) \cap C, N[u] \cap C\}$ and $\{N(v) \cap C, N[v] \cap C\}$ are disjoint for all distinct vertices $u,v \in V$.
  \\[-5mm]
\end{enumerate}
\end{theorem}

\begin{proof}
Let $u$ and $v$ be two distinct vertices of $G$.

\noindent (a) $\iff$ (b):
The result follows from the fact that, for any set $S \subset V$,
\begin{eqnarray}\label{eqn_PL}
\begin{split}
&(N(u) \cap S) \setminus \{v\} \ne (N(v) \cap S) \setminus \{u\}\\
\iff &(N(u) \cap S) \setminus \{u,v\} \ne (N(v) \cap S) \setminus \{u,v\}\\
\iff &\big( (N(u) \bigtriangleup N(v)) \setminus \{u,v\} \big) \cap S \ne \emptyset.
\end{split}
\end{eqnarray}

\noindent (b) $\iff$ (c): The result follows by Observation~\ref{apobs_nbhd_idty}.

\noindent (b) $\implies$ (d): The result again follows by Observation~\ref{apobs_nbhd_idty}.

\noindent (d) $\implies$ (b): First, let $u$ and $v$ be adjacent in $G$. Then, we have $(N[u] \bigtriangleup N[v]) \cap C \ne \emptyset$. Therefore, Observation~\ref{apobs_nbhd_idty} implies that $C$ has a non-empty intersection with $\big( N(u) \bigtriangleup N(v) \big) \setminus \{u,v\}$. We now assume that $u$ and $v$ are non-adjacent in $G$. Then, we have $(N(u) \bigtriangleup N(v)) \cap C \ne \emptyset$. Therefore, again Observation~\ref{apobs_nbhd_idty} implies that $C$ has a non-empty intersection with $\big( N(u) \bigtriangleup N(v) \big) \setminus \{u,v\}$. This proves the result.

\noindent (d) $\iff$ (e): The result follows by the definitions of closed- and open-separating sets.

\noindent (d) $\iff$ (f): The part (f) $\implies$ (d) follows immediately from the definitions. Let us therefore assume (d) and prove (f). In this case, using (d), it is only left to prove that $N(u) \cap C \ne N[v] \cap C$ for two distinct vertices $u,v \in V$. Let us first assume that $u \notin C$ (respectively, $v \notin C$). Then we have $N(u) \cap C = N[u] \cap C$ (respectively, $N(v) \cap C = N[v] \cap C$). Therefore, again by using (d), we have $N(u) \cap C \ne N[v] \cap C$ and hence, the result holds in this case. Let us now assume that $u,v \in C$. Moreover, we assume to the contrary that $N(u) \cap C = N[v] \cap C$. In other words, $N[v] \cap C \subseteq N[u] \cap C$. Moreover, this also implies that $uv$ is an edge in $G$ and hence, $u \in N[v]$. Therefore, we have $N[u] \cap C \subseteq N[v] \cap C$ as well, that is, $N[u] \cap C = N[v] \cap C$ contradicting the assumption (d). This proves the result.
\end{proof}

%%%%%%%%%%%%%%%%%

Following the characterizations of a full-separating set in Theorem~\ref{thm_PL-char}(d), for any vertex subset $C$ of a graph $G$ and any two distinct vertices $u, v$ of $G$, we say that the set $C$ \emph{full-separates} the pair $u,v$ if $(N(u) \bigtriangleup N(v)) \cap C \ne \emptyset$ and $(N[u] \bigtriangleup N[v]) \cap C \ne \emptyset$. Moreover, by Theorem~\ref{thm_PL-char}, the hyperedge sets of both $\hyp_\FD(G)$ and $\hyp_\FTD(G)$ are composed of $N[G]$ for domination and $N(G)$ for total-domination, respectively, as well as $\bigtriangleup_a[G]$ and $\bigtriangleup_n(G)$ in both cases. In addition, it turns out  that a full-separating set incorporates both elements of closed-separation and open-separation. Thus, Theorem~\ref{thm_PL-char} also justifies how the full-separation property models monitoring systems that tackle both types of errors in 3-state detectors described in the introduction, whereby such detectors can turn into 2-state detectors, either of Type~1 or of Type~2. Finally, Theorem~\ref{thm_PL-char}(f) shows that the notions of FTD-codes and strongly identifying codes\footnote{A \emph{strongly identifying code} (or \emph{SID-code} for short) was defined in~\cite{H_2010} to be a total-dominating set $S$ of a graph $G$ such that the sets $\{N(u) \cap S, N[u] \cap S\}$ and $\{N(v) \cap S, N[v] \cap S\}$ are distinct for all distinct vertices $u,v$ of $G$.} (SID-codes)~\cite{H_2010} are equivalent.

%%%%%%%%%%%%%%%%%%%%%%%%%%%%%%%%%%%%%%%%%%%%%

\subsection{Existence of FD- and FTD-codes}

Let $\codes$ = \{LD, LTD, ID, ITD, OD, OTD, FD, FTD\}. It is known from the literature that the already studied X-codes for some $\X \in \codes$ may not exist in all graphs, see for example \cite{CW_ISCO2024,HHH_2006,KCL_1998,SS_2010}. Considering the existence of X-codes from the point of view of covers of the X-hypergraphs $\hyp_\X(G) = (V,\F_\X)$, we immediately see that $G$ has no X-code if and only if $\hyp_\X(G)$ has no cover if and only if $\emptyset \in \F_\X$. In this context, it has been observed in \cite{ABLW_2018,ABLW_2022} that there is no X-code in a graph $G$ with $\hyp_\X(G)$ involving \\[-5mm]
%%%%%%%%%%%%%%%%%%%%%%
\begin{itemize}[leftmargin=12pt, itemsep=0pt]
\item $N(G)$ if $G$ has \emph{isolated vertices}, that is, vertices $v$ with $N(v) = \emptyset$;
\item $\bigtriangleup_a[G]$ if $G$ has \emph{closed twins}, that is, adjacent vertices $u,v$ with $N[u] = N[v]$ 
(and thus $N[u] \bigtriangleup N[v] = \emptyset$);
\item $\bigtriangleup_n(G)$ if $G$ has \emph{open twins}, that is, two non-adjacent vertices $u,v$ with $N(u) = N(v)$ 
(and thus $N(u) \bigtriangleup N(v) = \emptyset$). \\[-5mm]
\end{itemize}
%%%%%%%%%%%%%%%%%%%%%%
\begin{table}[ht]
\begin{center}
\begin{tabular}{ c || c | c || c | c || c | c || c | c }
 X & LD   & LTD  & ID   & ITD  & OD & OTD & FD & FTD \\
 \hline
 &&&&&&&&\\
 & $N[G]$ & \framebox[1.26cm]{$N(G)$} & $N[G]$ & \framebox[1.26cm]{$N(G)$}  & $N[G]$ & \framebox[1.26cm]{$N(G)$} & $N[G]$ & \framebox[1.26cm]{$N(G)$}\\
 &&&&&&&&\\
$\F_\X$ & $\bigtriangleup_a(G)$ & $\bigtriangleup_a(G)$ 
 & \framebox[1.26cm]{$\bigtriangleup_a[G]$}~ & \framebox[1.26cm]{$\bigtriangleup_a[G]$}  & $\bigtriangleup_a(G)$ & $\bigtriangleup_a(G)$ & \framebox[1.26cm]{$\bigtriangleup_a[G]$} & \framebox[1.26cm]{$\bigtriangleup_a[G]$}\\
 &&&&&&&&\\
 &  $\bigtriangleup_n[G]$ & $\bigtriangleup_n[G]$ & $\bigtriangleup_n[G]$ & $\bigtriangleup_n[G]$  & \framebox[1.26cm]{$\bigtriangleup_n(G)$} & \framebox[1.26cm]{$\bigtriangleup_n(G)$} & \framebox[1.26cm]{$\bigtriangleup_n(G)$} & \framebox[1.26cm]{$\bigtriangleup_n(G)$}\\
\end{tabular}
\end{center}
\caption{
The X-hypergraphs $\hyp_\X(G) = (V,\F_\X)$, listing in boxes the subsets of hyperedges which may contain an empty hyperedge and hence, may  
lead to non-existence of X-codes.}
\label{tab_admissible}
\end{table}
%%%%%%%%%%%%%%%%%%%%%%
Table \ref{tab_admissible} illustrates which X-problems for $\X \in \codes$ may not exist in a graph due to the presence of an empty set in the collections mentioned in the boxes. Calling a graph $G$ \emph{X-admis\-sible} if $G$ has an X-code, we see, for example, from Table \ref{tab_admissible} that every graph $G$ is LD-admissible, whereas a graph $G$ is OTD-admissible if and only if $G$ has neither isolated vertices nor open twins. Moreover, we call a graph \emph{twin-free} if it has neither closed nor open twins. Accordingly, we conclude the following regarding the existence of FD- and FTD-codes in graphs:
%%%%%%%%%%%%%%%%%%%%%%%%
\begin{corollary}\label{cor_admissible}
A graph $G$ is \\[-5mm]
\begin{itemize}[leftmargin=12pt, itemsep=0pt]
\item \FD-admissible if and only if $G$ is twin-free; 
\item \FTD-admissible if and only if $G$ is twin-free and has no isolated vertex. \\[-5mm]
\end{itemize}
\end{corollary}
%%%%%%%%%%%%%%%%%%%%%%%%
Since any two distinct isolated vertices of a graph are open twins with the empty set as both their open neighborhoods, Corollary~\ref{cor_admissible} further implies the following.
%%%%%%%%%%%%%%%%%%%%%%%%
\begin{corollary}\label{cor_FDadmis_1isol}
An \FD-admissible graph has at most one isolated vertex.
\end{corollary}
%%%%%%%%%%%%%%%%%%%%%%%%

Moreover, a special interest lies in hyperedges of $\hyp_\X(G) = (V,\F_\X)$ consisting of a single vertex, called an \emph{$\X$-forced vertex}, as each such vertex must belong to every \X-code of $G$. We denote the set of singleton hyperedges containing the \X-forced vertices of $G$ by $\for_{\X}(G)$. In \cite{ABLW_2018,ABLW_2022}, it was examined which of the six possible hyperedge subsets $\F'$ composing $\F_\X$ can contain singleton hyperedges. Denoting by $\for_{\X}(\F')$ the subset of singleton hyperedges of $\F'$, it has been observed in \cite{ABLW_2018,ABLW_2022} that for a graph $G=(V,E)$, 
\begin{itemize}
\item $\for(N[G])$ consists of all isolated vertices of $G$; 
\item $\for(N(G))$ contains all support vertices $v$ of pendant vertices $u$ of $G$, that is, all vertices $v$ such that there exists $u \in V \setminus \{v\}$ with $N(u)=\{v\}$; 
\item $\for(\bigtriangleup_a[G])$ is composed of all vertices $v \in V$ such that there exist adjacent vertices $u,u' \in V \setminus \{v\}$ with $N[u] \bigtriangleup N[u']=\{v\}$;
\item $\for(\bigtriangleup_n(G))$ consists of all vertices $v \in V$ such that there exist non-adjacent vertices $u,u' \in V \setminus \{v\}$ with $N(u) \bigtriangleup N(u')=\{v\}$;
  \end{itemize}
whereas $\for(\bigtriangleup_a(G))=\for(\bigtriangleup_n[G])=\emptyset$.
Accordingly, we say that all vertices of a graph $G$ belonging to
\begin{itemize}
\item $\for(N[G])$ are \emph{domination-forced}, 
\item $\for(N(G))$ are \emph{total-domination-forced}, 
\item $\for(\bigtriangleup_a[G]) \cup \for(\bigtriangleup_n(G))$ are \emph{full-separation-forced}, 
\end{itemize}
and conclude:

\begin{corollary}\label{cor_forced}
For a graph $G$, we have 
\[
\begin{array}{rcl}
\for_{\FD}(G) & = & \for(N[G]) \cup \for(\bigtriangleup_a[G]) \cup \for(\bigtriangleup_n(G)); \text{ and} \\
\for_{\FTD}(G) & = & \for(N(G)) \cup \for(\bigtriangleup_a[G]) \cup \for(\bigtriangleup_n(G)). \\
\end{array}
\]
\end{corollary}

In view of the hypergraph representations of the \X-codes with $\X \in \codes$, we see that $\for(\bigtriangleup_a[G])$ equals exactly the set of vertices forced by closed-separation, whereas $\for(\bigtriangleup_n(G))$ equals exactly the set of vertices forced by open-separation. Hence, the set of all full-separation-forced vertices also reflects the fact that full-separation combines both aspects of closed- and open-separation. 

\begin{example} \rm
For illustration, we construct $\hyp_{\FD}(P_4)$ and $\hyp_{\FTD}(P_4)$. 
The involved hyperedges are the following: 
\[
\begin{array}{lc|lc|lc|l}
N[G]           & & N(G)         & & \bigtriangleup_a[G]                & & \bigtriangleup_n(G) \\ \hline
N[1]=\{1,2\}   & & N(1)=\{2\}   & & N[1] \bigtriangleup N[2] = \{3\}   & & N(1) \bigtriangleup N(3) = \{4\}\\
N[2]=\{1,2,3\} & & N(2)=\{1,3\} & & N[2] \bigtriangleup N[3] = \{1,4\} & & N(1) \bigtriangleup N(4) = \{2,3\}\\
N[3]=\{2,3,4\} & & N(3)=\{2,4\} & & N[3] \bigtriangleup N[4] = \{2\}   & & N(2) \bigtriangleup N(4) = \{1\}\\
N[4]=\{3,4\}   & & N(4)=\{3\}   & &                            & & \\
\end{array}
\]
Accordingly, we see that $\for(N[P_4])=\emptyset$, $\for(N(P_4))=\bigtriangleup_a[P_4])=\{2,3\}$ and
$\for(\bigtriangleup_n(P_4)) = \{1,4\}$ holds.
This shows that $\for_{FD}(P_4) = \for_{FTD}(P_4) = \{1,2,3,4\}$ and thus $\fd(P_4) = \ftd(P_4) = 4$ follows.
\end{example}

Further, recall that each \fdadmis ~graph has at most one isolated vertex. 
Moreover, for an \fdadmis ~graph $G$ having an isolated vertex $v$, $\for(\bigtriangleup_n(G))$ contains all vertices $u$ of degree 1 (as $N(u) \bigtriangleup N(v) = N(u) \bigtriangleup \emptyset = N(u)$ holds).

In this context, we note that $\hyp_\X(G) = (V,\F_\X)$ may contain redundant hyperedges, see also \cite{ABLW_2018,ABLW_2022,ABW_2016,CW_ISCO2024}. 
In fact, if there are two hyperedges $F, F' \in \F_\X$ with $F \subseteq F'$, then $F \cap C \neq \emptyset$ also implies $F' \cap C \neq \emptyset$ for every $C \subseteq V$. 
Thus, $F'$ is \emph{redundant} as $(V,\F_\X -\{F'\})$ suffices to encode the X-codes of $G$. 
Hence, only non-redundant hyperedges of $\hyp_\X(G)$ need to be considered in order to determine $\tau(\hyp_\X(G))$ and thus $\x(G)$. 
In particular, if one hyperedge $F \in \F_\X$ contains an \X-forced vertex of $G$, then $F$ is clearly redundant.

\subsection{Relations of X-numbers and resulting bounds on FD- and FTD-numbers}

Consider the set $\mathcal{H}^*(V)$ of all hypergraphs $\hyp = (V,\mathcal{F})$ on the same vertex set $V$. We define a relation $\prec $ on $\mathcal{H}^*(V) \times \mathcal{H}^*(V)$ by $\hyp \prec \hyp'$ if, for every hyperedge $F'$ of $\hyp'$, there exists a hyperedge $F$ of $\hyp$ such that $F \subseteq F'$. We can show:
%%%%%%%%%%%%%%%%%%%%%%%
\begin{lemma}\label{lem_pos}
If $\hyp, \hyp' \in \mathcal{H}^*(V)$ with $\hyp \prec \hyp'$, then any cover of $\hyp$ is also a cover of $\hyp'$. In particular, we have $\tau(\hyp') \leq \tau(\hyp)$.
\end{lemma}
%%%%%%%%%%%%%%%%%%%%%%%
\begin{proof}
Let $\hyp, \hyp' \in \mathcal{H}^*(V)$ with $\hyp = (V,\mathcal{F})$, $\hyp' = (V,\mathcal{F}')$ and assume that $\hyp \prec \hyp'$. Moreover, let $C \subseteq V$ be a cover of $\hyp$. Now, let $F' \in \F'$ be any hyperedge. Since $\hyp \prec \hyp'$, there exists a hyperedge $F$ of $\hyp$ such that $F \subseteq F'$. Then $C \cap F \neq \emptyset$, as $C$ is a cover of $\hyp$. This implies that $C \cap F' \neq \emptyset$ as well. Hence, $C$ intersects every hyperedge of $\hyp'$ and thus, is also a cover of $\hyp'$. The second statement follows by the fact that, if $C$ is a minimum cover of $\hyp$, then $C$ is also a cover of $\hyp'$ and hence, $\tau(\hyp') \le |C| = \tau(\hyp)$.
\end{proof}
%%%%%%%%%%%%%%%%%%%%%%

\begin{figure}[!t]
\begin{center}
\includegraphics[scale=0.8]{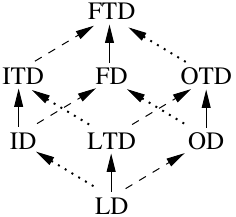}
\caption{The relations between the X-numbers for all $X \in \codes$, where $X' \longrightarrow X$ stands for $\gamma^{\X'}(G) \leq \gamma^\X(G)$ and solid (respectively, dotted and dashed) arrows refer to case (a) (respectively, (b) and (c)) of Corollary \ref{cor_all_relations}.}
\label{fig_relations}
\end{center}
\end{figure}

It is easy to see that for every graph $G = (V,E)$, the X-hypergraphs $\hyp_\X(G) = (V,\F_\X)$ belong to $\mathcal{H}^*(V)$ for all $X \in \codes$. Furthermore, it is clear that \\[-5mm]
%%%%%%%%%%%%%%%%%%%%%%
\begin{itemize}[leftmargin=12pt, itemsep=0pt]
\item $N[v] = N(v) \cup \{v\}$ for all vertices $v \in V$ (by definition);
\item $N(u) \bigtriangleup N(v) = (N[u] \bigtriangleup N[v]) \cup \{u,v\}$ for all adjacent vertices $u,v \in V$;
\item $N[u] \bigtriangleup N[v] = (N(u) \bigtriangleup N(v)) \cup \{u,v\}$ for all non-adjacent vertices $u,v \in V$. \\[-5mm]
\end{itemize}
%%%%%%%%%%%%%%%%%%%%%%%
This observation combined with Lemma~\ref{lem_pos} implies: 
%%%%%%%%%%%%%%%%%%%%%%%
\begin{corollary}\label{cor_all_relations}
If for a graph $G = (V,E)$ and two problems $\X,\X' \in \codes$, the hyperedges of $\hyp_\X(G) = (V,\F_\X)$ and $\hyp_{\X'}(G) = (V,\F_{\X'})$ \emph{only} differ in either the  \\[-5mm]
%%%%%%%%%%%%%%%%%%%%%%%
\begin{itemize}[leftmargin=18pt, itemsep=0pt]
\item[(a)] neighborhoods such that $N(G) \subset \F_\X$ and $N[G] \subset \F_{\X'}$, or
\item[(b)] the symmetric differences of neighborhoods of adjacent vertices such that $\bigtriangleup_a[G] \subset \F_\X$ and $\bigtriangleup_a(G) \subset \F_{\X'}$, or
\item[(c)] the symmetric differences of neighborhoods of non-adjacent vertices such that 
$\bigtriangleup_n(G) \subset \F_\X$ and $\bigtriangleup_n[G] \subset \F_{\X'}$, \\[-5mm]
\end{itemize}
then we have $\hyp_\X(G) \prec \hyp_{\X'}(G)$ in all the above cases which implies in particular $\gamma^{\X'}(G) \leq \gamma^\X(G)$.
\end{corollary}
%%%%%%%%%%%%%%%%%%%%%%%%%
From Corollary \ref{cor_all_relations} combined with the hypergraph representations from Table \ref{tab_admissible}, we immediately obtain the relations between the X-numbers for all $X \in \codes$ shown in Figure \ref{fig_relations}. 
%%%%%%%%%%%%%%%%%%%%%%%%
%%%%%%%%%%%%%%%%%%%%%%%%
Note that most of the relations not involving FD- or FTD-numbers were already known before, see, for example, \cite{ABLW_2022,CW_ISCO2024}. However, from the global point of view on all X-numbers given in Figure \ref{fig_relations}, we see in particular that $\ld(G)$ is a lower bound for \emph{all} other X-numbers, whereas $\ftd(G)$ is an upper bound for \emph{all} other X-numbers in FTD-admissible graphs.

As is evident from Figure~\ref{fig_relations}, the FD- and FTD-codes are, in general, stronger than, for example, ID-, OD-, LD-codes, etc. However, it is worth mentioning here that in the literature of identification problems, such notions of stronger codes related to stronger monitoring requirements in networks have already been presented by other authors. Some of the more prominent ones among them being the notions of \emph{strongly identifying codes}~\cite{H_2010} (which turn out to be equivalent to FTD-codes; see Theorem~\ref{thm_PL-char} for a proof of this equivalence), \emph{strongly $(t,\le l)$-identifying codes}~\cite{HLR_2002}, \emph{robust identifying codes} and \emph{dynamic identifying codes}~\cite{HKL_2006}, those of \emph{fault-tolerant} LD-, ID- and OTD-codes in three principle varieties, namely \emph{redundant}, \emph{error-detecting} and \emph{error-correcting} introduced, for example, in~\cite{JS_2022, JS_2022_ID1, JS_2022_ID2, JS_2024, JS_2025} and the notions of \emph{self locating-dominating sets} and \emph{solid locating-dominating sets}~\cite{JLL_2023}.

%%%%%%%%%%%%%%%%%%%%%%%%%%%%%%%%%%%%%%%%%%%%%%%

The following remark shows that in order to check if a total-dominating set $C$ of a graph $G$ is an \FTD-code of $G$, we do not need to check if $C$ full-separates every pair of distinct vertices of $G$ but only those which are at a distance of at most~$2$ between them. Here, the \emph{distance} between two vertices $u$ and $v$ of $G$, denoted by $d(u,v)$, is the length of a shortest path connecting $u$ and $v$ in $G$.

\begin{remark} \label{rem:FTD_dist_2}
Let $G$ be an \FTD-admissible graph. A total-dominating set $C$ of $G$ is an \FTD-code of $G$ if and only if $C$ full-separates every pair $u,v$ of distinct vertices of $G$ such that $d(u,v) \le 2$.
\end{remark}

\begin{proof}
The necessary condition for the statement follows immediately from the definition of an  \FTD-code. We, therefore, prove the sufficient condition. Thus, it is enough to show that $C$ full-separates every pair $u,v$ of distinct vertices of $G$ such that $d(u,v) \ge 3$. So, assume $u$ and $v$ to be a pair of distinct vertices of $G$ such that $d(u,v) \ge 3$. Since $C$ is a total-dominating set of $G$, the vertex $u$ has a neighbor, say $w$, in $C$. However, $w$ is not a neighbor of $v$ (since $d(u,v) \ge 3$) and hence, $w$ belongs to both $(N(u) \bigtriangleup N(v)) \cap C$ and $(N[u] \bigtriangleup N[v]) \cap C$. Thus, $C$ full-separates $u,v$ and hence, this proves the result.
\end{proof}

The following remark is immediate from the definition of full-separating sets.

\begin{remark} \label{rem:FD_empty_signature}
Let $G$ be an \FD-admissible graph and let $C$ be an \FD-code of $G$. Then, there exists at most one vertex $w \in V(G)$ such that $N(w) \cap C = \emptyset$.
\end{remark}

An equivalence of Remark~\ref{rem:FTD_dist_2} does not entirely apply for \FD-codes. In other words, it is not true that any dominating set of a graph $G$ that full-separates distinct vertices of $G$ of distance at most~$2$ is an \FD-code of $G$. As the next remark shows, for such a result to be true for \FD-codes, we need the additional condition that the dominating set is ``almost'' a total-dominating set, that is, it total-dominates all vertices of $G$ except possibly for one.

\begin{remark} \label{rem:FD_dist_2}
Let $G$ be an \FD-admissible graph and let $C$ be a dominating set of $G$ such that there exists at most one vertex $w \in V(G)$ with $N(w) \cap C = \emptyset$. If the set $C$ full-separates all distinct vertices $u,v \in V(G)$ with $d(u,v) \leq 2$, then $C$ is an \FD-code of $G$.
\end{remark}

\begin{proof}
It is enough to show that $C$ full-separates every pair $u,v$ of distinct vertices of $G$ such that $d(u,v) \geq 3$. In other words, since $uv \notin E$, it is enough to show that $(N(u) \bigtriangleup N(v)) \cap C \ne \emptyset$. Therefore, for two such vertices $u,v \in V$ with $d(u,v) \geq 3$, we have $N(u) \triangle N(v) = N(u) \cup N(v)$. If to the contrary, $(N(u) \bigtriangleup N(v)) \cap C = \emptyset$, then it implies that $(N(u) \cup N(v)) \cap C = \emptyset$ and so, $N(u) \cap C = N(v) \cap C = \emptyset$. However, this contradicts our assumption that $C$ is a full-separating set of $G$. This proves the result.
\end{proof}

The next theorem demonstrates the closeness between the
\FD- and the \FTD-number of a graph and hence, motivates us to compare the two numbers on graphs of some specific families later in Section \ref{sec4} to verify if they differ or are equal to each other on these graphs.

\begin{theorem} \label{thm:FTD=FD+1} \label{thm_pd-ptd_relations}
Let $G$ be an \FD-admissible graph. If $G$ is a disjoint union of a graph $G'$ and an isolated vertex, then we have
\[
\fd (G) = \ftd (G') + 1.
\] 
Otherwise, we have
\[
\ftd (G) - 1 \leq \fd (G) \leq \ftd(G).
\]
\end{theorem}

\begin{proof}
Let us first assume that $G$ is a disjoint union of a graph $G'$ and an isolated vertex $v$. Moreover, let $C$ be a minimum \FD-code of $G$. Then, $v \in C$ in order for $C$ to dominate $v$. This implies that $N(v) \cap C = \emptyset$. Therefore, by Remark~\ref{rem:FD_dist_2}, the set $C' = C \setminus \{v\}$ must total-dominate all vertices in $V \setminus \{v\}$. Moreover, as $v$ is an isolated vertex of $G$, the set $C'$ is also a full-separating set of $G'$ and hence, an \FTD-code of $G'$. This implies that $\ftd(G') \leq |C| - 1 = \fd(G)-1$, that is, $\fd(G) \geq \ftd(G') + 1$. On the other hand, if $C'$ is a minimum \FTD-code of $G'$, then $C' \cup \{v\}$ is an \FD-code of $G$. This implies that $\fd(G) \leq |C'|+1 = \ftd(G')+1$ and hence, the result follows.

\medskip

Let us now assume that the graph $G$ has no isolated vertices. Then, $G$ is also \FTD-admissible. Now, the inequality $\fd (G) \leq \ftd (G)$ holds since every \FTD-code is also an \FD-code. To prove the other inequality, we let $C$ be a minimum \FD-code of $G$. Then, we have $|C| = \fd(G)$. Now, if $C$ is a total-dominating set of $G$, then $C$ is also an \FTD-code of $G$ and hence, the result holds trivially. Therefore, let us assume that $C$ is not a total-dominating set of $G$. This implies that there exists a vertex $v$ such that $N(v) \cap C = \emptyset$. Moreover, Remark~\ref{rem:FD_empty_signature} implies that $v$ is the only vertex of $G$ whose neighborhood has an empty intersection with $C$. Since $G$ has no isolated vertices by assumption, the vertex $v$ has a neighbor, say $u$, in $G$. Then, $C \cup \{u\}$ is a total-dominating set of $G$ and hence, is also a full-separating set of $G$. Therefore, we have $\ftd (G) \leq |C|+1 = \fd (G) +1$.
\end{proof}

%%%%%%%%%%%%%%%%%%%%%%%%%%%%%%%%%%%%%%%%%%%%%%%

\noindent
It is clear from Figure~\ref{fig_relations} that for any \FTD-admissible graph $G$, we have $\ftd(G) \ge \x(G)$, where $\X \in \{\ITD, \OTD, \FD\}$ and for any \FD-admissible graph $G$, we have $\fd(G) \ge \x(G)$, where $\X \in \{\ID, \OD\}$. Finally, using Theorem~\ref{thm:FTD=FD+1} and the relations in Figure~\ref{fig_relations}, we have $\fd(G) \ge \ftd(G) - 1 \ge \x(G) - 1$, where $\X \in \{\LTD, \ITD, \OTD\}$. This leads us to the following result on lower bounds for FD- and FTD-numbers of graphs in terms of other \X-numbers for $\X \in \codes$.
%%%%%%%%%%%%%%%%%%%%%%%%
\begin{corollary}\label{cor_lb_PL-numbers}
Let $G$ be an \FTD-admissible graph. We have \\[-2mm]
\[
\begin{array}{lcl}
\ftd (G) & \geq & \max(\itd (G),\otd (G),\fd (G)), \\[1mm]
\fd (G)  & \geq & \max(\id (G),\od (G),\ltd (G)-1,\itd (G)-1,\otd (G)-1). \\
\end{array}
\]
\end{corollary}

%%%%%%%%%%%%%%%%%%%%%%%%

\section{Complexity of the FD- and FTD-problems} \label{sec3} 

In this section, we address the computational complexities of the \FD- and \FTD-problems. The decision versions of these two problems are formally stated as follows.

\defproblem{\mincode{\FD}}{$(G,k)$: A graph $G$ and a positive integer $k$.}{Does there exist an FD-code $S$ of $G$ such that $|S| \leq k$?}

\defproblem{\mincode{\FTD}}{$(G,k)$: A graph $G$ and a positive integer $k$.}{Does there exist an FTD-code $S$ of $G$ such that $|S| \leq k$?}

%%%%%%%%%%%%%%%%%%%%%%%%
We prove the above two problems to be NP-complete. Moreover, despite a difference of at most one as established in Theorem~\ref{thm_pd-ptd_relations}, we show that it is hard to decide if the FD- and FTD-numbers of a given input graph are the same or different. The decision version of this problem is defined as follows.

\defproblem{\FD ~$\ne$ \FTD}{A graph $G$ and an integer $k$.}{Is $\ftd(G) = k$ and $\fd(G) = k-1?$ That is, are the following assertions true?
\begin{enumerate}[leftmargin=18pt, itemsep=0pt]
\item There exists an \FTD-code $S$ of $G$ such that  $|S| = k$.
\item There exists an \FD-code $S'$ of $G$ such that $|S'| = k-1$.
\item For any vertex subset $S''$ of $G$, if $|S''| < |S|$, then $S''$ is not an \FTD-code of $G$; and if $|S''| < |S'|$, then $S''$ is not an \FD-code of $G$.
\end{enumerate}
}

As can be noticed in \FD~$\ne$ \FTD, given a vertex subset $S$ (with $|S| = k$) of an input graph $G$ on $n$ vertices, 
to check condition (3), one has to consider all subsets of $V(G)$ of order $k-1$ (which admits $\calO(n^{k-1})$ as a worst-case running time) and check if $S$ is an \FTD-code (which admits $\calO(n^2)$ as a running time). In other words, verifying a certificate for \FD~$\ne$ \FTD\ can potentially take up to a running time of order $\calO(n^{k-1}) \cdot \calO(n^2) \subseteq \calO(n^{k+1})$ which may not be polynomial in $n+k$. This implies that, \FD~$\ne$ \FTD\ does not necessarily belong to the class \NP. However, as we shall see in this section, the problem still remains \NP-hard.

%%%%%%%%%%%%%%%%%%%%%%%%

To prove the hardness results for the above three problems, we use one single reduction from an instance $\psi$ of 3-SAT to an instance $(G^\psi,k)$ of \mincode{\FD}, \mincode{\FTD} and FD $\ne$ FTD. 3-SAT happens to be the famous \NP-complete problem proven to be so by Karp in~\cite{K_1972}. For all standard terminologies and notations in the theory of computational complexity, we refer to the book~\cite{GJ_1990} by Garey and Johnson. From here on, we shall denote by $\psi = (X, \calC)$ a generic instance of 3-SAT  defined on the set $X$ of Boolean variables and the set $\calC$ of clauses defined on the literals from $X$. In addition, let $n = |X|$ and $m = |\calC|$ and let us denote by $\phi = \phi(\calC)$ the Boolean formula defined by $\phi = \mathbf{c_1} \wedge \mathbf{c_2} \wedge \ldots \wedge \mathbf{c_m}$, where $\mathbf{c_i} \in \calC$ for all $1 \le i \le m$. We may also assume from here on that no clause contains both literals of a variable, as such a clause is always satisfiable. 

%%%%%%%%%%%%%%%%%%%%%%%%%%%%%%%%%%%%%%%%%%%%%%%%%%%%%%%%%

\begin{reduction} \label{red_NP}
The reduction takes as input an instance $\psi = (X,\calC)$ of \textsc{3-SAT}. Let $n = |X|$ and $m = |\calC|$. The reduction constructs a graph $G^\psi$ on $10n+3m$ vertices as follows (refer to Figure~\ref{fig_NP}):
%%%%%%%%%%%%%%%%%%%%%%%%
\begin{itemize}[leftmargin=12pt, itemsep=0pt]
\item For every variable $x \in X$, do the following (refer to Figure~\ref{fig_NP}(a)):
%%%%%%%%%%%%%%%%%%%%%%%
\begin{itemize}[leftmargin=8pt, itemsep=0pt]
\item Add 10 vertices named $v^x_1$, $v^x_2$, $v^x_3$, $w^x_1$, $w^x_2$, $s^x_1$, $s^x_2$, $s^x_3$, $z^x_1$ and $z^x_2$ (intuitively, the vertex $w^x_1$ will correspond to the literal $x$ and the vertex $w^x_2$ will correspond to the literal $\neg x$).
\item Add edges so that the graph induced by the vertices $v^x_3$, $w^x_1$, $w^x_2$, $s^x_1$, $s^x_2$ and $s^x_3$ is a $K_6$ minus the edges $w^x_1w^x_2$, $v^x_1s^x_1$, $v^x_1s^x_2$ and $v^x_1s^x_3$.
\item Add edges $v^x_1v^x_2$ and $v^x_2v^x_3$ making the vertices $v^x_1$, $v^x_2$ and $v^x_3$ induce a $P_3$; and add edges $s^x_1 z^x_1$ and $s^x_2z^x_2$ making $z^x_1$ and $z^x_2$ pendant vertices and the vertices $s^x_1$ and $s^x_2$ their respective supports.
\end{itemize}
%%%%%%%%%%%%%%%%%%%%%%%%
Let the graph induced by the above 10 vertices be denoted as $G^x$ and be called the \emph{variable gadget} corresponding to the variable $x \in X$.
%%%%%%%%
\item For every clause $\mathbf{c} \in \calC$, do the following (refer to Figure~\ref{fig_NP}(b)):
%%%%%%%%%%%%%%%%%%%%%%%%%
\begin{itemize}[leftmargin=8pt, itemsep=0pt]
\item Add 3 vertices named $u^{\mathbf{c}}_1$, $u^{\mathbf{c}}_2$ and $u^{\mathbf{c}}_3$.
\item Add edges $u^{\mathbf{c}}_1u^{\mathbf{c}}_2$ and $u^{\mathbf{c}}_2u^{\mathbf{c}}_3$ thus making $u^{\mathbf{c}}_1$, $u^{\mathbf{c}}_2$ and $u^{\mathbf{c}}_3$ induce a $P_3$.
\end{itemize}
%%%%%%%%%%%%%%%%%%%%%%%%%
Let the $P_3$ induced by the above 3 vertices be denoted by $P^\mathbf{c}$ and be called the \emph{clause gadget} corresponding to the clause $\mathbf{c} \in \calC$.
%%%%%%%%%
\item For all variables $x \in X$ and all clauses $\mathbf{c} \in \calC$, if the literal $x$ is in a clause $\mathbf{c}$, then add the edge $u^\mathbf{c}_1w^x_1$; and if the literal $\neg x$ is in a clause $\mathbf{c}$, then add the edge $u^\mathbf{c}_1w^x_2$ (refer to Figure~\ref{fig_NP}(c)).
\end{itemize}
\end{reduction}
%%%%%%%%%%%%%%%%%%%%%%%%%%%
\begin{figure}[!t]
\centering
\begin{subfigure}[t]{0.25\textwidth}
\centering
\begin{tikzpicture}[scale=0.33,
blacknode/.style={circle, draw=black!, fill=black!, thick, minimum size= 2mm, scale=0.6},
whitenode/.style={circle, draw=black!, fill=white!, thick, minimum size= 2mm, scale=0.6},
]

\node[whitenode] (1) at (0,0) {}; \node at (-0.7,0.7) {$w^x_1$};
\node[whitenode] (2) at (4,0) {}; \node at (4.7,0.7) {$w^x_2$};
\node[whitenode] (3) at (2,2) {}; \node at (2.7,2.7) {$v^x_1$};
\node[whitenode] (4) at (2,4) {}; \node at (2.7,4.7) {$v^x_2$};
\node[whitenode] (5) at (2,6) {}; \node at (2.7,6.7) {$v^x_3$};
\node[whitenode] (6) at (0,-2) {}; \node at (-0.7,-2.7) {$s^x_1$};
\node[whitenode] (7) at (4,-2) {}; \node at (4.7,-2.7) {$s^x_2$};
\node[whitenode] (8) at (2,-4) {}; \node at (2,-5) {$s^x_3$};
\node[whitenode] (9) at (0,-4) {}; \node at (-0.7,-4.7) {$z^x_1$};
\node[whitenode] (10) at (4,-4) {}; \node at (4.7,-4.7) {$z^x_2$};

\draw[-, thick, black!] (1) -- (3);
\draw[-, thick, black!] (1) -- (6);
\draw[-, thick, black!] (1) -- (7);
\draw[-, thick, black!] (1) -- (8);
\draw[-, thick, black!] (2) -- (3);
\draw[-, thick, black!] (2) -- (6);
\draw[-, thick, black!] (2) -- (7);
\draw[-, thick, black!] (2) -- (8);
\draw[-, thick, black!] (3) -- (4);
\draw[-, thick, black!] (4) -- (5);
\draw[-, thick, black!] (6) -- (7);
\draw[-, thick, black!] (6) -- (8);
\draw[-, thick, black!] (7) -- (8);
\draw[-, thick, black!] (6) -- (9);
\draw[-, thick, black!] (7) -- (10);

\end{tikzpicture}
\caption{$G^x$: Variable gadget corresponding to $x \in X$.} \label{fig_NP Variable}
\end{subfigure}
\hspace{2mm}
\begin{subfigure}[t]{0.25\textwidth}
\centering
\begin{tikzpicture}[scale=0.33,
blacknode/.style={circle, draw=black!, fill=black!, thick, minimum size= 2mm, scale=0.6},
whitenode/.style={circle, draw=black!, fill=white!, thick, minimum size= 2mm, scale=0.6},
]
\node[whitenode] (1) at (0,4) {}; \node at (1,4) {$u^\mathbf{c}_1$};
\node[whitenode] (2) at (0,2) {}; \node at (1,2) {$u^\mathbf{c}_2$};
\node[whitenode] (3) at (0,0) {}; \node at (1,0) {$u^\mathbf{c}_3$};

\draw[-, thick, black!] (1) -- (2);
\draw[-, thick, black!] (2) -- (3);

\end{tikzpicture}
\caption{$P^\mathbf{c}$: Clause gadget corresponding to $\mathbf{c} \in \calC$.} \label{fig_NP Clause}
\end{subfigure}
\hspace{2mm}
\begin{subfigure}[t]{0.35\textwidth}
\centering
\begin{tikzpicture}[scale=0.33,
blacknode/.style={circle, draw=black!, fill=black!, thick, minimum size= 2mm, scale=0.6},
whitenode/.style={circle, draw=black!, fill=white!, thick, minimum size= 2mm, scale=0.6},
]

\node[blacknode] (1) at (0,0) {}; \node at (-0.7,0.7) {$w^x_1$};
\node[whitenode] (2) at (4,0) {}; \node at (4.7,0.7) {$w^x_2$};
\node[blacknode] (3) at (2,2) {}; \node at (2.7,2.7) {$v^x_1$};
\node[blacknode] (4) at (2,4) {}; \node at (2.7,4.7) {$v^x_2$};
\node[whitenode] (5) at (2,6) {}; \node at (2.7,6.7) {$v^x_3$};
\node[blacknode, fill=black!30] (6) at (0,-2) {}; \node at (-0.7,-2.7) {$s^x_1$};
\node[blacknode] (7) at (4,-2) {}; \node at (4.7,-2.7) {$s^x_2$};
\node[whitenode] (8) at (2,-4) {}; \node at (2,-5) {$s^x_3$};
\node[blacknode] (9) at (0,-4) {}; \node at (-0.7,-4.7) {$z^x_1$};
\node[blacknode] (10) at (4,-4) {}; \node at (4.7,-4.7) {$z^x_2$};

\node[blacknode] (11) at (-3,-4) {}; \node at (-4,-4) {$u^\mathbf{c}_1$};
\node[blacknode] (12) at (-3,-6) {}; \node at (-4,-6) {$u^\mathbf{c}_2$};
\node[whitenode] (13) at (-3,-8) {}; \node at (-4,-8) {$u^\mathbf{c}_3$};

\draw[-, thick, black!] (1) -- (3);
\draw[-, thick, black!] (1) -- (6);
\draw[-, thick, black!] (1) -- (7);
\draw[-, thick, black!] (1) -- (8);
\draw[-, thick, black!] (2) -- (3);
\draw[-, thick, black!] (2) -- (6);
\draw[-, thick, black!] (2) -- (7);
\draw[-, thick, black!] (2) -- (8);
\draw[-, thick, black!] (3) -- (4);
\draw[-, thick, black!] (4) -- (5);
\draw[-, thick, black!] (6) -- (7);
\draw[-, thick, black!] (6) -- (8);
\draw[-, thick, black!] (7) -- (8);
\draw[-, thick, black!] (6) -- (9);
\draw[-, thick, black!] (7) -- (10);

\draw[-, thick, black!] (11) -- (12);
\draw[-, thick, black!] (12) -- (13);

\draw[-, thick, black!] (11) -- (1);
\draw[dashed, thick, black!] (11) -- (-4,-0.5);
\draw[dashed, thick, black!] (11) -- (-6,-0.5);
\draw[dashed, thick, black!] (1) -- (-1.7,-1.5);

\draw[dashed, thick, black!] (2) -- (6,-1.5);

\end{tikzpicture}
\caption{$G^\psi$: Instance of \mincode{\FD}, \mincode{\FTD} and \FD~$\ne$ \FTD.} \label{fig_NP G_I}
\end{subfigure}
\caption[Reduction~\ref{red_NP} from \textsc{3-SAT} to \mincode{\FD}, \mincode{\FTD} and \FD~$\ne$ \FTD]{Polynomial-time construction of the graph $G^\psi$ from an instance $\psi = (X,\calC)$ of \textsc{3-SAT} as in Reduction~\ref{red_NP}. The black vertices in (c) represent those in a code described in Lemma~\ref{aplem_NP_1}. The gray vertex ($s^x_1$) implies that, for some fixed variable $x=x' \in X$, the vertex is not included in an \FD-code but is included in the \FTD-code described in Lemma~\ref{aplem_NP_1}.}
\label{fig_NP}
\end{figure}
%%%%%%%%%%%%%%%%%%%%%%%%
It can be verified that Reduction~\ref{red_NP} is carried out in time polynomial in $m$ and $n$.

\begin{lemma} \label{aplem_NP}
Let $S$ be a full-separating set of $G^\psi$.

\begin{enumerate}[leftmargin=18pt, itemsep=0pt]
\item \label{aplem_NP_Py} If $S$ is either an \FD- or an \FTD-code of $G^\psi$, then for all $\mathbf{c} \in \calC$, we have $|V(P^\mathbf{c}) \cap S| \ge 2$.

\item \label{aplem_NP_Gx_FTD} If $S$ is an \FTD-code, then we have $|(V(G^x) \setminus \{w^x_1,w^x_2\}) \cap S| \ge 6$ for all $x \in X$.

\item \label{aplem_NP_Gx_FD} If $S$ is an \FD-code, then except possibly for one $x' \in X$ for which $|( V(G^{x'}) \setminus \{w^{x'}_1,w^{x'}_2\}) \cap S| \ge 5$, we have $|(V(G^x) \setminus \{w^x_1,w^x_2\}) \cap S| \ge 6$ for all $x \in X \setminus \{x'\}$.
\end{enumerate}

Moreover, we have $\ftd(G^\psi) \ge 7n+2m$ and $\fd(G^\psi) \ge 7n+2m-1$.
\end{lemma}

\begin{proof}
(1): Let $S$ be either an \FD- or an \FTD-code of $G^\psi$ and let $\mathbf{c} \in \calC$. Since $u^\mathbf{c}_3$ is a pendant vertex of $G^\psi$, we must have $|\{u^\mathbf{c}_2,u^\mathbf{c}_3\} \cap S| \ge 1$. Moreover, the vertex $u^{\mathbf{c}}_1$ is full-separation-forced with respect to the pair $u^{\mathbf{c}}_2,u^{\mathbf{c}}_3$. Therefore, we also have $u^{\mathbf{c}}_1 \in S$. This proves (1).

\medskip

(2) and (3): Let $S$ be either an \FD- or an \FTD-code of $G^\psi$ and let $x \in X$. Since $v^x_3$ is also a pendant vertex of $G^\psi$, we must have $|\{v^x_2,v^x_3\} \cap S| \ge 1$. Since the vertex $v^x_1$ is full-separation-forced with respect to the pair $v^x_2,v^x_3$ and, for $i \in \{1,2\}$ the vertex $z^x_i$ is full-separation-forced with respect to $s^x_i,s^x_3$, we must have $v^x_1,z^x_1,z^x_2 \in S$ for all $x \in X$. Furthermore, we have the following.
\begin{itemize}[leftmargin=12pt, itemsep=0pt]
\item If $S$ is an \FTD-code, for each $i \in \{1,2\}$ and each $x \in X$, the vertex $s^x_i$ is a support vertex of $G^\psi$ and hence, is total-domination-forced. Therefore, we have $s^x_i \in S$.
\item If $S$ is an \FD-code, by Remark~\ref{rem:FD_empty_signature}, at least all but one, that is, $(n-1)$ many, vertices in the set $\{s^x_i \in V(G^\psi): x \in X, i \in \{1,2\}\}$ 
must belong to $S$.
\end{itemize}
This proves parts (2) and (3).

\medskip

Finally, for all $x \in X$, in order for $S$ to full-separate the pair $v_1^x, v_3^x$, we must have $|\{w_1^x, w_2^x\} \cap S| \ge 1$. Hence, by counting the vertices in $S$, the final statement of the result follows.
\end{proof}

\begin{lemma}\label{aplem_NP_1}
If $X$ has a satisfying assignment, then we have $\ftd(G^\psi) = 7n+2m$ and $\fd(G^\psi) = 7n+2m-1$.
\end{lemma}

\begin{proof}
First, we show that there exists an \FTD-code $S$ of $G^\psi$ such that $|S| = 7n+2m$. To that end, we build a vertex subset $S$ of $G^\psi$ as follows.
\begin{itemize}[leftmargin=12pt, itemsep=0pt]
\item For all $x \in X$ and $\mathbf{c} \in \calC$, pick the vertices $u^\mathbf{c}_1, u^\mathbf{c}_2,v^x_1,v^x_2,z^x_1,z^x_2,s^x_1, s^x_2$ in $S$.
\item For any variable $x \in X$, either $x$ or $\neg x$ is assigned~$1$ in the satisfying assignment on $X$. If $x$ is assigned~$1$, then pick $w^x_1$ in $S$; and if $\neg x$ is assigned~$1$, then pick $w^x_2$ in $S$.
\end{itemize}
Note that, by construction, we have $|S| = 7n+2m$. We now show that $S$ is an \FTD-code of $G^\psi$. Notice that, again by construction, the set $S$ is a total-dominating set of $G^\psi$. Therefore, by Remark~\ref{rem:FTD_dist_2}, it is enough to show that $S$ full-separates all pairs of distinct vertices of $G^\psi$ with distance at most 2 between them. For every $x \in X$ and $\mathbf{c} \in \calC$, the set $S$ full-separates the vertices
\begin{itemize}[leftmargin=12pt, itemsep=0pt]
\item $v^x_2, v^x_3$ by $v^x_1$; the vertices $u^\mathbf{c}_2, u^\mathbf{c}_3$ by $u^\mathbf{c}_1$; and the vertices $v^x_1, u^\mathbf{c}_1$ by $u^\mathbf{c}_2$;
\item $v^x_1, v^x_3$, by either $w^x_1$ or $w^x_2$ according to whether $(x, \neg x) = (1,0)$ or $(x, \neg x) = (0,1)$, respectively;
\item $u^\mathbf{c}_1, u^\mathbf{c}_3$ by either $w^x_1$ or $w^{x'}_2$ for some $x,x' \in X$ according to whether $x$ (assigned to 1) belongs to the clause $\mathbf{c}$; or $\neg x'$ (assigned to 1) belongs to $\mathbf{c}$;
\item $w^x_1, w^x_2$ by $u^\mathbf{c}_1$ for some clause $\mathbf{c} \in \calC$ (to which exactly one of $x$ and $\neg x$ belongs);
\item $w^x_i, s^x_j$ by $v^x_1$ for any $i \in \{1,2\}$ and $j \in \{1,2,3\}$;
\item $w^x_i, z^x_j$ by $v^x_1$; the vertices $w^x_i,u^\mathbf{c}_j$ by $v^x_1$; and the vertices $w^x_i,v^\mathbf{c}_j$ by $s^x_2$ for any $i,j \in \{1,2\}$;
\item $s^x_i, s^x_j$, by either $z^x_1$ or $z^x_2$ for any $i,j \in \{1,2,3\}$;
\item $s^x_i, v^x_1$, by $z^x_i$; and the vertices $s^x_i, u^\mathbf{c}_1$, by $z^x_i$ for any $i \in \{1,2,3\}$;
\item $s^x_i, z^x_j$, by either $w^x_1$ or $w^x_2$ according to whether $(x, \neg x) = (1,0)$ or $(x, \neg x) = (0,1)$, respectively, where $i,j \in \{1,2\}$.
\end{itemize}
Finally, the only vertices in different variable gadgets (respectively, clause gadgets) with distance 2 between them are of the form $w^x_i \in G^x$ and $w^{x'}_j \in G^{x'}$ (respectively, $u^\mathbf{c}_1 \in P^\mathbf{c}$ and $u^{\mathbf{c}'}_1 \in P^{\mathbf{c}'}$), where $x,x' \in X$ (respectively, $\mathbf{c},\mathbf{c}' \in \calC$) are distinct. However, the set $S$ full-separates any such pair $w^x_i, w^{x'}_j$ by $v^x_1$ and the pair $u^\mathbf{c}_1, u^{\mathbf{c}'}_1$ by $u^\mathbf{c}_2$. This proves that $S$ is a full-separating set of $G^\psi$.

\medskip

We now let $S' = S \setminus \{s^{x'}_1\}$, where $x'$ is a fixed variable in $X$. Therefore, we have $|S'| = 7n+2m-1$. Hence, it is enough to show that $S'$ is an \FD-code of $G^\psi$. First of all, it can be verified that $S'$ is a dominating set of $G^\psi$. We now show that $S'$ is also a full-separating set of $G^\psi$. To do so, notice that $N_{G^\psi}(z^{x'}_1) \cap S' = \emptyset$ for only the variable $x' \in X$. Therefore, by Remark~\ref{rem:FD_dist_2}, it is enough to show that $S'$ full-separates all pairs of distinct vertices of $G^\psi$ with distance at most 2 between them. In other words, we consider for full-separation by $S'$ the exact same pairs of vertices we had considered for full-separation by $S$ in the previous list. Moreover, we notice that in the above list showing how $S$ full-separates pairs of vertices of $G^\psi$, we see that no two vertices are uniquely separated by $s^x_1 \in S$ for any $x \in X$. This implies that each of those pairs of vertices of $G^\psi$ is full-separated by $S'$ as well and hence, this proves that $S'$ too is a full-separating set of $G^\psi$.

\medskip

Thus, the existence of the above \FTD-code $S$ and the \FD-code $S'$ implies that $\ftd(G^\psi) \le |S| = 7n+2m$ and $\fd(G^\psi) \le |S'| = 7n+2m-1$. Therefore, combining these with the last part of Lemma~\ref{aplem_NP}, the result follows.
\end{proof}

\begin{lemma} \label{aplem_NP_w_sat_assign}
Let $\psi = (X,\calC)$ be an instance of \textsc{$3$-SAT} and $G^\psi$ be the graph as in Reduction~\ref{red_NP}. If there exists a full-separating set $S$ of $G^\psi$ such that $|\{w^x_1,w^x_2\} \cap S| = 1$ for all $x \in X$, then $X$ has a satisfying assignment.
\end{lemma}

\begin{proof}
Let $S$ be a full-separating set of $G^\psi$ such that $|\{w^x_1,w^x_2\} \cap S| = 1$ for all $x \in X$. We now provide a Boolean assignment on $X$ the following way: for any $x \in X$, if $w^x_1 \in S$ and $w^x_2 \notin S$, then put $(x,\neg x)= (1,0)$; and if $w^x_1 \notin S$ and $w^x_2 \in S$, then put $(x,\neg x) = (0,1)$. Clearly, this assignment on $X$ is a valid one since, for each $x \in X$, exactly one of $x$ and $\neg x$ is assigned $1$ and the other $0$. To now prove that this assignment on $X$ is also a satisfying one, we simply note that, for each clause $\mathbf{c} \in \calC$, in order for the set $S$ to full-separate $u^\mathbf{c}_1$ and $u^\mathbf{c}_3$, either there exists a variable $x \in X$ such that $u^\mathbf{c}_1w^x_1 \in E(G^\psi)$ with $w^x_1 \in S$; or there exists a variable $x' \in X$ such that $u^\mathbf{c}_1w^{x'}_2 \in E(G^\psi)$ and $w^{x'}_2 \in S$. Therefore, by construction of the graph $G^\psi$, either $x$ is a variable in the clause $\mathbf{c}$ with $(x,\neg x) = (1,0)$ or $x'$ is a variable whose literal $\neg x'$ is in the clause $\mathbf{c}$ with $(x', \neg x') = (0,1)$. Hence, the binary assignment formulated on $X$ is a satisfying one.
\end{proof}

\begin{lemma}\label{aplem_NP_2}
The following statements are true. 
\begin{enumerate}[leftmargin=18pt, itemsep=0pt]
\item\label{aplem_NP_2_FTD} If $\ftd(G^\psi) = 7n+2m$, then $X$ has a satisfying assignment.
\item\label{aplem_NP_2_FD} If $\fd(G^\psi) = 7n+2m-1$, then $X$ has a satisfying assignment.
\end{enumerate}
\end{lemma}

\begin{proof}
Let us first assume that $\ftd(G^\psi) = 7n+2m$.
Then, there exists an \FTD-code $S$ of $G^\psi$ such that $|S| = 7n+2m$. Then, we show that $X$ has a satisfying assignment. Therefore, by Lemma~\ref{aplem_NP_w_sat_assign}, we only need to show that $|\{w^x_1,w^x_2\} \cap S| = 1$ for all $x \in X$ (note that $|\{w^x_1,w^x_2\} \cap S| \ge 1$ in order for $S$ to full-separate the vertices $v^x_1,v^x_3$). If, on the contrary, $|\{w^{x^\star}_1, w^{x^\star}_2\} \cap S| = 2$ for some $x^\star \in X$, then, by Lemma~\ref{aplem_NP}(\ref{aplem_NP_Py})
and (\ref{aplem_NP_Gx_FTD}), this implies that
\begin{align*}
7n+2m = |S| &= \sum_{x \in X} |V(G^x) \cap S| + \sum_{\mathbf{c} \in B} |V(P^\mathbf{c}) \cap S|\\
&= \sum_{x \in X} |(V(G^x) \setminus \{w^x_1,w^x_2\}) \cap S| + \sum_{\mathbf{c} \in B} |V(P^\mathbf{c}) \cap S| + \sum_{x \in X} |\{w^x_1,w^x_2\} \cap S|\\
&\ge  6n+2m + |\{w^{x^\star}_1,w^{x^\star}_2\} \cap S| + \sum_{x \in X \setminus \{x^\star\}} |\{w^x_1,w^x_2\} \cap S| \ge 7n+2m+1.
\end{align*}
However, this is a contradiction. Therefore, we must have $|\{w^x_1,w^x_2\} \cap S| = 1$ for all $x \in X$. This proves (1).

\medskip

We now assume that $\fd(G^\psi) = 7n+2m-1$. Therefore, there exists an \FD-code $S'$ of $G^\psi$ such that $|S'| = 7n+2m-1$. To show that $X$ has a satisfying assignment, again by Lemma~\ref{aplem_NP_w_sat_assign}, we only need to show that $|\{w^x_1,w^x_2\} \cap S'| = 1$ for all $x \in X$. On the contrary, let us assume that $|\{w^{x^\star}_1, w^{x^\star}_2\} \cap S'| = 2$ for some $x^\star \in X$. Moreover, for all $x \in X$, we have $|\{w^x_1,w^x_2\} \cap S'| \ge 1$ in order for $S'$ to full-separate the vertices $v^x_1,v^x_3$. Therefore, using Lemma~\ref{aplem_NP}(\ref{aplem_NP_Gx_FD}), we must have $|V(G^x) \cap S'| \ge 7$ for all $x \in X$. This implies that
\begin{align*}
7n+2m-1 = |S'| &= \sum_{x \in X} |V(G^x) \cap S'| + \sum_{\mathbf{c} \in B} |V(P^\mathbf{c}) \cap S'| \ge 7n+2m.
\end{align*}
This is again a contradiction. Therefore, we must have $|\{w^x_1,w^x_2\} \cap S'| = 1$ for all $x \in X$ and this proves (2).
\end{proof}

This brings us to the proofs of the main results in this section.

\medskip
\begin{theorem} \label{thm_NP_FTD}
\mincode{\FTD} is \NP-complete.
\end{theorem}

\begin{proof}
\mincode{\FTD} clearly belongs to the class \NP\ since it can be verified in polynomial-time if a given vertex subset of a graph $G$ is an \FTD-code of $G$. To prove \NP-hardness, we show that an instance $\psi=(X,\calC)$ with $|X|=n$ and $|\calC|=m$ is a \yes-instance of \textsc{3-SAT} if and only if $(G^\psi, 7n+2m)$ is a \yes-instance of \mincode{\FTD}, where the graph $G^\psi$ is as constructed in Reduction~\ref{red_NP}. In other words, we show that $X$ has a satisfying assignment if and only if there exists an \FTD-code $S$ of $G^\psi$ such that $|S| \leq 7n+2m$.

\medskip

To prove the necessary part of the last statement, if $X$ has a satisfying assignment, then by Lemma~\ref{aplem_NP_1}, we have $\ftd(G^\psi) = 7n+2m$. This implies that there exists an \FTD-code $S$ of $G^\psi$ such that $|S| = 7n+2m$. On the other hand, to prove the sufficiency, if there exists an \FTD-code $S$ of $G^\psi$ such that $|S| \le 7n+2m$, then combining with Lemma~\ref{aplem_NP}, we have $\ftd(G^\psi) = 7n+2m$. Therefore, by Lemma~\ref{aplem_NP_2}(1), there exists a satisfying assignment on $X$. This proves the result.
\end{proof}
%%%%%%%%%%%%%%%%%%%%%%%%%%

%\medskip
\begin{theorem} \label{thm_NP_FD}
\mincode{\FD} is \NP-complete.
\end{theorem}

\begin{proof}
\mincode{\FD} clearly belongs to the class \NP\ since it can be verified in polynomial-time if a given vertex subset of a graph $G$ is an \FD-code of $G$. To prove \NP-hardness, we show that an instance $\psi=(X,\calC)$ with $|X|=n$ and $|\calC|=m$ is a \yes-instance of \textsc{3-SAT} if and only if $(G^\psi,7n+2m-1)$ is a \yes-instance of \mincode{\FD}, where the graph $G^\psi$ is as constructed in Reduction~\ref{red_NP}. In other words, we show that $X$ has a satisfying assignment if and only if there exists an \FD-code $S$ of $G^\psi$ such that $|S| \leq 7n+2m-1$.

\medskip

To prove the necessary part of the last statement, if $X$ has a satisfying assignment, then by Lemma~\ref{aplem_NP_1}, we have $\fd(G^\psi) = 7n+2m-1$. This implies that there exists an \FD-code $S$ of $G^\psi$ such that $|S| = 7n+2m-1$. On the other hand, to prove the sufficiency, if there exists an \FD-code $S$ of $G^\psi$ such that $|S| \le 7n+2m-1$, then combining with Lemma~\ref{aplem_NP}, we have $\fd(G^\psi) = 7n+2m-1$. Therefore, by Lemma~\ref{aplem_NP_2}(2), there exists a satisfying assignment on $X$. This proves the result.
\end{proof}
%%%%%%%%%%%%%%%%%%%%%%%%%%

\begin{theorem} \label{thm_FD neq FTD}
\FD ~$\ne$ \FTD ~is \NP-hard.
\end{theorem}

\begin{proof}
We show that \FD ~$\ne$ \FTD\ is \NP-hard by showing that an instance $\psi=(X,\calC)$ with $|X|=n$ and $|\calC|=m$ is a \yes-instance of \textsc{3-SAT} if and only if $(G^\psi, 7n+2m)$ is a \yes-instance of FD $\ne$ FTD, where the graph $G^\psi$ is as constructed in Reduction~\ref{red_NP}. In other words, we show that $X$ has a satisfying assignment if and only if 
$\fd(G^\psi) = 7n+2m-1$ and $\ftd(G^\psi) = 7n+2m$.

\medskip

The necessary condition is true since, if $X$ has a satisfying assignment, then by Lemma~\ref{aplem_NP_1}, we have $\ftd(G^\psi) = 7n+2m$ and $\fd(G^\psi) = 7m+2m-1$. To prove the sufficiency condition, let us assume that $\ftd(G^\psi) = 7n+2m$ and $\fd(G^\psi) = 7n+2m-1$. Then, by Lemma~\ref{aplem_NP_2}, $X$ has a satisfying assignment.
\end{proof}

%%%%%%%%%%%%%%%%%%%%%%%%%%%%%%%%%%%%%%

\section{\FD- and \FTD-numbers of some graph families} \label{sec4}

In this section, we study the \FD- and the \FTD-numbers of graphs belonging to some well-known graph families. Moreover, motivated by Theorem~\ref{thm:FTD=FD+1} showing that the \FD-number and the \FTD-number of a graph differ by at most~$1$, we compare the two numbers on graphs of these families.

\subsection{Paths and cycles}
 
In order to study the \FD- and \FTD-numbers of paths $P_n$ and cycles $C_n$, we first note that $P_2$ and $C_3$ have closed twins; and $P_3$ and $C_4$ have open twins. However, apart from them, all other paths and cycles are \FD- and \FTD-admissible. The proof of the following theorem giving the exact values of the \FD- and \FTD-numbers for paths and cycles is based on results for the closed formulas of the \OTD-numbers of paths given in~\cite{SS_2010} and of cycles given in~\cite{BCLW_2024}, respectively.

\begin{lemma}\label{lem_pd=ptd paths}
Let $G$ be either a path $P_n$ for $n \geq 4$ or a cycle $C_n$ for $n \geq 5$. Then, we have $\fd(G) = \ftd(G)$.
\end{lemma}

\begin{proof}
The graph $G$ has maximum degree~$2$. On the contrary, let us assume that $\fd(G) \neq \ftd(G)$. Moreover, let $C$ be a minimum \FD-code of $G$. Therefore, we have $|C| = \fd(G)$. Therefore, $C$ is not a total-dominating set of $G$ (or else $C$ would be an \FTD-code of $G$ thus implying $\fd(G) = \ftd(G)$, a contradiction). Now, by Remark~\ref{rem:FD_dist_2}, there exists exactly one vertex $w$ of $G$ such that $N(w) \cap C = \emptyset$. Thus, $w \in C$ in order for $C$ to dominate $w$. Let us first assume that $w$ is not a full-separation-forced vertex. Since $N(w) \cap C = \emptyset$, every vertex $w' \in N(w)$ has a neighbor other than $w$ in $C$ (or else, $w$ and $w'$ are not full-separated by $C$, a contradiction). Therefore, for any neighbor $w'$ of $w$ in $G$, the set $C' = (C \setminus \{w\}) \cup \{w'\}$ is a total-dominating set of $G$. Moreover, since $w \in C$ is not full-separation-forced, the set $C'$ is also a full-separating set of $G$ and hence, is an \FTD-code of $G$. Moreover, we have $|C'| = |C| = \fd(G)$. This implies that $\ftd(G) \le \fd(G)$ and hence, $\fd(G) = \ftd(G)$, contrary to our assumption.

\medskip

Therefore, let us assume that the vertex $w$ is full-separation-forced with respect to $u$ and $v$. This implies, without loss of generality, that $u \in N(w)$ and $v \notin N[w]$. Moreover, since $N(w) \cap C = \emptyset$ and $v \neq w$, by Remark~\ref{rem:FD_empty_signature}, we must have $N(v) \cap C \neq \emptyset$. Therefore, let $x \in N(v) \cap C$. Then, $x \notin N[w]$. Now, since $C$ full-separates $u$ and $v$ uniquely by $w$ (since $w$ is full-separation-forced), we must have $x \in N(u)$ as well. Since $u \notin C$, in order for $C$ to full-separate $v$ and $x$, there must exist another vertex $y \in (N(v) \setminus N(x)) \cap C$ (notice that $y$ cannot be a neighbor of $x$ since $\{u,v\} \subset N(x)$ and $\deg(x) \le 2$). Moreover, $uy \notin E$ since $\{w,x\} \subset N(u)$ and $\deg(u) \le 2$. Thus, the set $C$ full-separates $u$ and $v$ by $y$, a contradiction to the uniqueness of $w$. Hence, we must have $\fd(G) = \ftd(G)$.
\end{proof}

\begin{observation} \label{obs_S geq 4}
Let $G=(V,E)$ be either a path $P_n$ for $n \geq 4$ or a cycle $C_n$ for $n \geq 5$. Let $C$ be an \FTD-code of $G$. Then, for any vertex $v \in V$, there exists a vertex subset $S_v \subset C$ of $G$ such that $v \in S_v$, $|S_v| \geq 4$ and $G[S_v]$ is the path $P_{|S_v|}$.
\end{observation}

\begin{proof}
Let $V = \{v_1, v_2, \ldots , v_n\}$ such that $v_iv_{i+1} \in E$ for all $i \in \{1, \ldots , n-1\}$ and $v_n v_1 \in E$ if $G \cong C_n$. Let $v=v_i \in C$ for some $i \in \{1, \ldots , n\}$. In order for $C$ to total-dominate $v_i$, we must have $v_{i+1} \in C$, say, without loss of generality. Next, in order for $C$ to full-separate $v_i$ and $v_{i+1}$, we must have $v_{i+2} \in C$, say, without loss of generality. Again, in order for $C$ to full-separate $v_i$ and $v_{i+2}$, we must have $v_{i+3} \in C$, say, without loss of generality. Thus, the result follows by taking $S = \{v_i,v_{i+1},v_{i+2},v_{i+3}\}$.
\end{proof}
\begin{lemma} \label{lem_pd_ptd_paths-cycles_ub}
Let $G$ be either a path $P_n$ for $n \geq 4$ or a cycle $C_n$ for $n \geq 5$. Let $n=6q+r$, where $q$ and $r$ are non-negative integers with $r \in \{1,2, \ldots , 5\}$.
Then we have
\[
\ftd (G) \leq \begin{cases}4q+r, &\text{if $r \in \{1, \ldots , 4\}$;}\\
4q+4, &\text{if $r=5$}. \end{cases}
\]
\end{lemma}

\begin{proof}
Let $V = V(G) = \{v_1, v_2, \ldots , v_n\}$ and assume that $n = 6q+r$, where $r \in \{0,1,2,3,4,5\}$. If $q=0$, the graph $G$ is either a  $4$-path or a $5$-path or a $5$-cycle. In all the three cases, it can be checked that the set $\{v_1, v_2, v_3, v_4\}$ is an \FTD-code of $G$. Thus, in these cases, the result holds. For the rest of this proof, therefore, we assume that $n \geq 6$, that is, $q \geq 1$. We now construct a vertex subset $C$ of $G$ by including in $C$ the vertices $v_{6k-4}, v_{6k-3}, v_{6k-2}, v_{6k-1}$ for all $k \in \{1,2, \ldots , q\}$, the vertices $v_{6q}, v_{6q+1}, \ldots , v_{6q+r-1}$ if $r \in \{1,2,3,4\}$ and the vertices $v_{6q+1}, v_{6q+2}, v_{6q+3}, v_{6q+4}$ if $r=5$ (see Figure~\ref{fig_paths} for an example with the code vertices marked in black). Notice that, by construction, the vertex subset $C$ is a total-dominating set of $G$. Therefore, by Remark~\ref{rem:FTD_dist_2}, to show that $C$ is an \FTD-code of $G$, it is enough to show that $C$ full-separates every pair of distinct $u,v \in V$ such that $d(u,v) \leq 2$. Now, one can observe that for every pair $v_i, v_j$ of distinct vertices in $G$ such that $i<j \leq i+2$ (that is $d (v_i,v_j) \leq 2$), either $v_{i-1}$ or $v_{j+1}$ (or both) is in $C$. Thus, all vertex pairs $v_i, v_j \in V$ with $d(v_i,v_j) \leq 2$ are full-separated by $C$. This implies that $C$ is an \FTD-code of $G$. Hence, by counting the order of $C$ from its construction, the result follows by $\ftd(G) \leq |C|$.
\end{proof}
%
%

%%%%%%%%%%%%%%%%%%%%%%%%
\begin{figure}[t]
\centering
\begin{tikzpicture}[scale=0.3,
blacknode/.style={circle, draw=black!, fill=black!, thick, minimum size= 2mm, scale=0.5},
whitenode/.style={circle, draw=black!, fill=white!, thick, minimum size= 2mm, scale=0.5},
]

\node[whitenode] (1) at (0,0) {}; \node at (0,0.8) {$v_1$};
\node[blacknode] (2) at (2,0) {}; \node at (2,0.8) {$v_2$};
\node[blacknode] (3) at (4,0) {}; \node at (4,0.8) {$v_3$};
\node[blacknode] (4) at (6,0) {}; \node at (6,0.8) {$v_4$};
\node[blacknode] (5) at (8,0) {}; \node at (8,0.8) {$v_5$};
\node[whitenode] (6) at (10,0) {}; \node at (10,0.8) {$v_6$};
\node[whitenode] (7) at (12,0) {}; \node at (12,0.8) {$v_7$};
\node[blacknode] (8) at (14,0) {}; \node at (14,0.8) {$v_8$};
\node[blacknode] (9) at (16,0) {}; \node at (16,0.8) {$v_9$};
\node[blacknode] (10) at (18,0) {}; \node at (18,0.8) {$v_{10}$};
\node[blacknode] (11) at (20,0) {}; \node at (20,0.8) {$v_{11}$};
\node[whitenode] (12) at (22,0) {}; \node at (22,0.8) {$v_{12}$};

\node at (36,0) {$P_{12}$: $r=0$};

%%%%%%%%%%%%%%%%%%%%%%%%%%%%%%%%%%

\node[whitenode] (1') at (0,-3) {}; \node at (0,-2.2) {$v_1$};
\node[blacknode] (2') at (2,-3) {}; \node at (2,-2.2) {$v_2$};
\node[blacknode] (3') at (4,-3) {}; \node at (4,-2.2) {$v_3$};
\node[blacknode] (4') at (6,-3) {}; \node at (6,-2.2) {$v_4$};
\node[blacknode] (5') at (8,-3) {}; \node at (8,-2.2) {$v_5$};
\node[whitenode] (6') at (10,-3) {}; \node at (10,-2.2) {$v_6$};
\node[whitenode] (7') at (12,-3) {}; \node at (12,-2.2) {$v_7$};
\node[blacknode] (8') at (14,-3) {}; \node at (14,-2.2) {$v_8$};
\node[blacknode] (9') at (16,-3) {}; \node at (16,-2.2) {$v_9$};
\node[blacknode] (10') at (18,-3) {}; \node at (18,-2.2) {$v_{10}$};

\node[blacknode] (11') at (20,-3) {}; \node at (20,-2.2) {$v_{11}$};
\node[blacknode] (12') at (22,-3) {}; \node at (22,-2.2) {$v_{12}$};
\node[whitenode] (13') at (24,-3) {}; \node at (24,-2.2) {$v_{13}$};

\node at (36,-3) {$P_{13}$: $r=1$};

%%%%%%%%%%%%%%%%%%%%%%%%%%%%%%%%%%%

\node[whitenode] (1'') at (0,-6) {}; \node at (0,-5.2) {$v_1$};
\node[blacknode] (2'') at (2,-6) {}; \node at (2,-5.2) {$v_2$};
\node[blacknode] (3'') at (4,-6) {}; \node at (4,-5.2) {$v_3$};
\node[blacknode] (4'') at (6,-6) {}; \node at (6,-5.2) {$v_4$};
\node[blacknode] (5'') at (8,-6) {}; \node at (8,-5.2) {$v_5$};
\node[whitenode] (6'') at (10,-6) {}; \node at (10,-5.2) {$v_6$};
\node[whitenode] (7'') at (12,-6) {}; \node at (12,-5.2) {$v_7$};
\node[blacknode] (8'') at (14,-6) {}; \node at (14,-5.2) {$v_8$};
\node[blacknode] (9'') at (16,-6) {}; \node at (16,-5.2) {$v_9$};
\node[blacknode] (10'') at (18,-6) {}; \node at (18,-5.2) {$v_{10}$};

\node[blacknode] (11'') at (20,-6) {}; \node at (20,-5.2) {$v_{11}$};
\node[blacknode] (12'') at (22,-6) {}; \node at (22,-5.2) {$v_{12}$};
\node[blacknode] (13'') at (24,-6) {}; \node at (24,-5.2) {$v_{13}$};
\node[whitenode] (14'') at (26,-6) {}; \node at (26,-5.2) {$v_{14}$};

\node at (36,-6) {$P_{14}$: $r=2$};

%%%%%%%%%%%%%%%%%%%%%%%%%%%%%%%%%%%%%

\node[whitenode] (1''') at (0,-9) {}; \node at (0,-8.2) {$v_1$};
\node[blacknode] (2''') at (2,-9) {}; \node at (2,-8.2) {$v_2$};
\node[blacknode] (3''') at (4,-9) {}; \node at (4,-8.2) {$v_3$};
\node[blacknode] (4''') at (6,-9) {}; \node at (6,-8.2) {$v_4$};
\node[blacknode] (5''') at (8,-9) {}; \node at (8,-8.2) {$v_5$};
\node[whitenode] (6''') at (10,-9) {}; \node at (10,-8.2) {$v_6$};
\node[whitenode] (7''') at (12,-9) {}; \node at (12,-8.2) {$v_7$};
\node[blacknode] (8''') at (14,-9) {}; \node at (14,-8.2) {$v_8$};
\node[blacknode] (9''') at (16,-9) {}; \node at (16,-8.2) {$v_9$};
\node[blacknode] (10''') at (18,-9) {}; \node at (18,-8.2) {$v_{10}$};

\node[blacknode] (11''') at (20,-9) {}; \node at (20,-8.2) {$v_{11}$};
\node[blacknode] (12''') at (22,-9) {}; \node at (22,-8.2) {$v_{12}$};
\node[blacknode] (13''') at (24,-9) {}; \node at (24,-8.2) {$v_{13}$};
\node[blacknode] (14''') at (26,-9) {}; \node at (26,-8.2) {$v_{14}$};
\node[whitenode] (15''') at (28,-9) {}; \node at (28,-8.2) {$v_{15}$};

\node at (36,-9) {$P_{15}$: $r=3$};

%%%%%%%%%%%%%%%%%%%%%%%%%%%%%%%%%%%%%%

\node[whitenode] (1*) at (0,-12) {}; \node at (0,-11.2) {$v_1$};
\node[blacknode] (2*) at (2,-12) {}; \node at (2,-11.2) {$v_2$};
\node[blacknode] (3*) at (4,-12) {}; \node at (4,-11.2) {$v_3$};
\node[blacknode] (4*) at (6,-12) {}; \node at (6,-11.2) {$v_4$};
\node[blacknode] (5*) at (8,-12) {}; \node at (8,-11.2) {$v_5$};
\node[whitenode] (6*) at (10,-12) {}; \node at (10,-11.2) {$v_6$};
\node[whitenode] (7*) at (12,-12) {}; \node at (12,-11.2) {$v_7$};
\node[blacknode] (8*) at (14,-12) {}; \node at (14,-11.2) {$v_8$};
\node[blacknode] (9*) at (16,-12) {}; \node at (16,-11.2) {$v_9$};
\node[blacknode] (10*) at (18,-12) {}; \node at (18,-11.2) {$v_{10}$};

\node[blacknode] (11*) at (20,-12) {}; \node at (20,-11.2) {$v_{11}$};
\node[blacknode] (12*) at (22,-12) {}; \node at (22,-11.2) {$v_{12}$};
\node[blacknode] (13*) at (24,-12) {}; \node at (24,-11.2) {$v_{13}$};
\node[blacknode] (14*) at (26,-12) {}; \node at (26,-11.2) {$v_{14}$};
\node[blacknode] (15*) at (28,-12) {}; \node at (28,-11.2) {$v_{15}$};
\node[whitenode] (16*) at (30,-12) {}; \node at (30,-11.2) {$v_{16}$};

\node at (36,-12) {$P_{16}$: $r=4$};

%%%%%%%%%%%%%%%%%%%%%%%%%%%%%%%%%%%%%%%

\node[whitenode] (1**) at (0,-15) {}; \node at (0,-14.2) {$v_1$};
\node[blacknode] (2**) at (2,-15) {}; \node at (2,-14.2) {$v_2$};
\node[blacknode] (3**) at (4,-15) {}; \node at (4,-14.2) {$v_3$};
\node[blacknode] (4**) at (6,-15) {}; \node at (6,-14.2) {$v_4$};
\node[blacknode] (5**) at (8,-15) {}; \node at (8,-14.2) {$v_5$};
\node[whitenode] (6**) at (10,-15) {}; \node at (10,-14.2) {$v_6$};
\node[whitenode] (7**) at (12,-15) {}; \node at (12,-14.2) {$v_7$};
\node[blacknode] (8**) at (14,-15) {}; \node at (14,-14.2) {$v_8$};
\node[blacknode] (9**) at (16,-15) {}; \node at (16,-14.2) {$v_9$};
\node[blacknode] (10**) at (18,-15) {}; \node at (18,-14.2) {$v_{10}$};

\node[blacknode] (11**) at (20,-15) {}; \node at (20,-14.2) {$v_{11}$};
\node[whitenode] (12**) at (22,-15) {}; \node at (22,-14.2) {$v_{12}$};
\node[blacknode] (13**) at (24,-15) {}; \node at (24,-14.2) {$v_{13}$};
\node[blacknode] (14**) at (26,-15) {}; \node at (26,-14.2) {$v_{14}$};
\node[blacknode] (15**) at (28,-15) {}; \node at (28,-14.2) {$v_{15}$};
\node[blacknode] (16**) at (30,-15) {}; \node at (30,-14.2) {$v_{16}$};
\node[whitenode] (17**) at (32,-15) {}; \node at (32,-14.2) {$v_{17}$};

\node at (36,-15) {$P_{17}$: $r=5$};

%Lines

\draw[-, thick, black!] (1) -- (2);
\draw[-, thick, black!] (2) -- (3);
\draw[-, thick, black!] (3) -- (4);
\draw[-, thick, black!] (4) -- (5);
\draw[-, thick, black!] (5) -- (6);
\draw[-, thick, black!] (6) -- (7);
\draw[-, thick, black!] (7) -- (8);
\draw[-, thick, black!] (8) -- (9);
\draw[-, thick, black!] (9) -- (10);
\draw[-, thick, black!] (10) -- (11);
\draw[-, thick, black!] (11) -- (12);

%%%%%%%%%%%%%%%%%%%%%%%%%%%%%%%%%

\draw[-, thick, black!] (1') -- (2');
\draw[-, thick, black!] (2') -- (3');
\draw[-, thick, black!] (3') -- (4');
\draw[-, thick, black!] (4') -- (5');
\draw[-, thick, black!] (5') -- (6');
\draw[-, thick, black!] (6') -- (7');
\draw[-, thick, black!] (7') -- (8');
\draw[-, thick, black!] (8') -- (9');
\draw[-, thick, black!] (9') -- (10');
\draw[-, thick, black!] (10') -- (11');

\draw[-, thick, black!] (11') -- (12');
\draw[-, thick, black!] (12') -- (13');

%%%%%%%%%%%%%%%%%%%%%%%%%%%%%%%%%%%

\draw[-, thick, black!] (1'') -- (2'');
\draw[-, thick, black!] (2'') -- (3'');
\draw[-, thick, black!] (3'') -- (4'');
\draw[-, thick, black!] (4'') -- (5'');
\draw[-, thick, black!] (5'') -- (6'');
\draw[-, thick, black!] (6'') -- (7'');
\draw[-, thick, black!] (7'') -- (8'');
\draw[-, thick, black!] (8'') -- (9'');
\draw[-, thick, black!] (9'') -- (10'');
\draw[-, thick, black!] (10'') -- (11'');

\draw[-, thick, black!] (11'') -- (12'');
\draw[-, thick, black!] (12'') -- (13'');
\draw[-, thick, black!] (13'') -- (14'');

%%%%%%%%%%%%%%%%%%%%%%%%%%%%%%%%%%%

\draw[-, thick, black!] (1''') -- (2''');
\draw[-, thick, black!] (2''') -- (3''');
\draw[-, thick, black!] (3''') -- (4''');
\draw[-, thick, black!] (4''') -- (5''');
\draw[-, thick, black!] (5''') -- (6''');
\draw[-, thick, black!] (6''') -- (7''');
\draw[-, thick, black!] (7''') -- (8''');
\draw[-, thick, black!] (8''') -- (9''');
\draw[-, thick, black!] (9''') -- (10''');
\draw[-, thick, black!] (10''') -- (11''');

\draw[-, thick, black!] (11''') -- (12''');
\draw[-, thick, black!] (12''') -- (13''');
\draw[-, thick, black!] (13''') -- (14''');
\draw[-, thick, black!] (14''') -- (15''');

%%%%%%%%%%%%%%%%%%%%%%%%%%%%%%%%%%%%%

\draw[-, thick, black!] (1*) -- (2*);
\draw[-, thick, black!] (2*) -- (3*);
\draw[-, thick, black!] (3*) -- (4*);
\draw[-, thick, black!] (4*) -- (5*);
\draw[-, thick, black!] (5*) -- (6*);
\draw[-, thick, black!] (6*) -- (7*);
\draw[-, thick, black!] (7*) -- (8*);
\draw[-, thick, black!] (8*) -- (9*);
\draw[-, thick, black!] (9*) -- (10*);
\draw[-, thick, black!] (10*) -- (11*);

\draw[-, thick, black!] (11*) -- (12*);
\draw[-, thick, black!] (12*) -- (13*);
\draw[-, thick, black!] (13*) -- (14*);
\draw[-, thick, black!] (14*) -- (15*);
\draw[-, thick, black!] (15*) -- (16*);

%%%%%%%%%%%%%%%%%%%%%%%%%%%%%%%%%%%%%%

\draw[-, thick, black!] (1**) -- (2**);
\draw[-, thick, black!] (2**) -- (3**);
\draw[-, thick, black!] (3**) -- (4**);
\draw[-, thick, black!] (4**) -- (5**);
\draw[-, thick, black!] (5**) -- (6**);
\draw[-, thick, black!] (6**) -- (7**);
\draw[-, thick, black!] (7**) -- (8**);
\draw[-, thick, black!] (8**) -- (9**);
\draw[-, thick, black!] (9**) -- (10**);
\draw[-, thick, black!] (10**) -- (11**);

\draw[-, thick, black!] (11**) -- (12**);
\draw[-, thick, black!] (12**) -- (13**);
\draw[-, thick, black!] (13**) -- (14**);
\draw[-, thick, black!] (14**) -- (15**);
\draw[-, thick, black!] (15**) -- (16**);
\draw[-, thick, black!] (16**) -- (17**);

%%%%%%%%%%%%%%%%%%%%%%%%%%%%%%%%%%%%%%

\draw[dotted, thick, black!] (11,1) -- (11,-11);
\draw[dotted, thick, black!] (23,1) -- (23,-11);

\end{tikzpicture}
\caption[Examples of minimum \FTD-codes of paths $P_n$ with $n \ge 4$ and cycles $C_n$ with $n \ge 5$]{The set of black vertices represents an \FTD-code of a path $P_n$ with $n \ge 4$ and a cycle $C_n$ with $n \ge 5$ (by joining the vertices $v_1$ and $v_n$ by an edge in each case).} \label{fig_paths}
\end{figure}
%%%%%%%%%%%%%%%%%%%%%%%

\begin{theorem} \label{thm:FD-FTD_paths_cycles}
Let $G$ be either a path $P_n$ for $n \geq 4$ or a cycle $C_n$ for $n \geq 5$. Moreover, let $n=6q+r$ for non-negative integers $q$ and $r \in \{1,2, \ldots , 5\}$.
Then, we have
\[
\fd(G) = \ftd(G) = \begin{cases}4q+r, &\text{if $r \in \{1,2,3,4\}$};\\
4q+4, &\text{if $r=5$}.
\end{cases}
\]
\end{theorem}

\begin{proof}
Let $V = V(G) = \{v_1, v_2, \ldots , v_n\}$ such that $v_iv_{i+1} \in E = E(G)$ for all $i \in \{1,2, \ldots , n-1\}$ and $v_nv_1 \in E$ if $G$ is a cycle $C_n$. By Lemma~\ref{lem_pd=ptd paths}, it is enough to prove the result only for $\ftd(G)$. To begin with, let us assume that $G$ is a path $P_n$ for $n \geq 4$. Then, by Corollary~\ref{cor_lb_PL-numbers}, we have
\[
\ftd(P_n) \geq \otd(P_n) = \begin{cases}4q+r, &\text{if $r \in \{1,2,3,4\}$};\\
4q+4, &\text{if $r=5$}. \end{cases}
\]
The last inequality has been shown in~\cite{SS_2010}. Therefore, by Lemma~\ref{lem_pd_ptd_paths-cycles_ub}, the result holds for the \FTD-numbers of paths on $n \geq 4$ vertices.

\medskip

We now prove the result for $\ftd(G)$ when $G$ is a cycle $C_n$ for $n \geq 5$. By the result in~\cite{BCLW_2024}), we have
\begin{align} \label{eq_old}
\otd(C_n) = \begin{cases}4q+r, &\text{if $r \in \{0,1,2,4\}$}\\
4q+2, &\text{if $r=3$}\\
4q+4, &\text{if $r=5$}.
\end{cases}
\end{align}
Then again, by Lemma~\ref{lem_pd_ptd_paths-cycles_ub}, the result for $\ftd(C_n)$ holds for all $n = 6q+r$ except for $r=3$. Therefore, we next analyze only the $(6q+3)$-cycles. So, let $n = 6q+3$ for any positive integer $q$. Using Lemma~\ref{lem_pd_ptd_paths-cycles_ub}, we notice that it is enough to prove that $\ftd(C_n) \geq 4q+3$. So, let us assume on the contrary that there exists an \FTD-code $C$ of $G$ such that $|C| \leq 4q+2$.
\begin{claim}
For some $w \in V$, there exists a vertex subset $S_w \subset C$ of $G$ (as in Observation~\ref{obs_S geq 4}) such that $w \in S_u$, $|S_w| \geq 5$ and $G[S_w]$ is the path $P_{|S_w|}$.
\end{claim}

\begin{claimproof}
Without loss of generality, let us assume that, for every $v \in V$, the sets $S_v$ as in Observation~\ref{obs_S geq 4} are maximal. Contrary to the claim, let us assume that $|S_v| \leq 4$ for all $v \in V$. Then, again using Observation~\ref{obs_S geq 4}, we have $|S_v| = 4$, for all $v \in V$. Define $u \sim v$ if $S_u = S_v$. For some $i < j$ and $v_i \not \sim v_j$, if we have $\max \{i' : v_{i'} \in S_{v_i}\} = \min \{j' : v_{j'} \in S_{v_j}\}$, then it implies that $|S_{v_i}| \geq 8$, a contradiction. Hence, for every $i < j$ and $v_i \not \sim v_j$, there must be a vertex $v_k \in V \setminus C$ such that $\max \{i' : v_{i'} \in S_{v_i}\} < k < \min \{j' : v_{j'} \in S_{v_j}\}$. This implies that $|C| = 4q'$, where $q' \in [q]$. Now, for any $i < j$ and $v_i \not \sim v_j$, let $v_k, v_{k+1}, \ldots , v_{k+k'} \in V \setminus C$ such that $\max \{i' : v_{i'} \in S_{v_i}\} < k < k+1 < \ldots < k+k' < \min \{j' : v_{j'} \in S_{v_j}\}$. Then we must have $k' \in \{0,1\}$, or else, the vertex $v_{k+1}$ is not total-dominated by $C$, a contradiction. Therefore, by counting, we must have $n \leq 6q' \leq 6q$ which is a contradiction since $n=6q+3$. This proves our claim.
\end{claimproof} 

Now, let $S=S_w$ and $v_i, v_{i+1} \in S$. Then, let $G'$ be the $(6q+2)$-cycle obtained by contracting the edge $v_iv_{i+1}$. Let $V' = V(G') = \{u_1, u_2, \ldots , u_{6q+2}\}$, where $u_j = v_j$ for $j \in \{1,2, \ldots , i\}$ and $u_j = v_{j+1}$ for $j \in \{i+1, i+2, \ldots , 6q+2\}$. Moreover, let $C' = \{u_j : v_j \in C\}$. Then $C'$ is also a total-dominating set of $G'$. Moreover, we have $|S_{u_j}| \geq 4$ for all $j \in \{1,2, \ldots , 6q+2\}$. For any two vertices $u_j,u_{j'} \in V'$, where $j < j'$ and $j'-j \leq 2$, we can observe that either $u_{j-1}$ or $u_{j'+1}$ (or both) is in $C'$. Hence, by Remark~\ref{rem:FTD_dist_2}, the set $C'$ is an \FTD-code of $G'$. However, this is a contradiction to Equation (\ref{eq_old}), since now we have $\ftd(C_{6q+2}) \leq |C'| = 4q+1$. Hence, $|C| \geq 4q+3$ and this proves the result.
\end{proof}

%%%%%%%%%%%%%%%%%%%%%%%%%%

\subsection{Half-graphs} 

\begin{figure}[h]
\begin{center}
\includegraphics[scale=0.9]{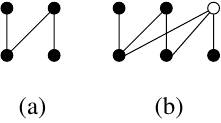}
\caption{Minimum FD-codes in half-graphs (the black vertices belong to the code), where (a) is $B_2=P_4$ and (b) is $B_3$.}
\label{fig_half-graph}
\end{center}
\end{figure}

A graph $G=(U \cup W, E)$ is \emph{bipartite} if its vertex set can be partitioned into two stable sets $U$ and $W$ so that every edge of $G$ has one endpoint in $U$ and the other in $W$. For any integer $k \geq 1$, the \emph{half-graph} $B_k=(U \cup W, E)$ is the bipartite graph with vertices in $U = \{u_1, \ldots , u_k\}$, $W = \{w_1, \ldots , w_k \}$ and edges $u_iw_j$ if and only if $i \leq j$ (see Figure~\ref{fig_half-graph}). In particular, we have $B_1=P_2$ and $B_2=P_4$. Moreover, we clearly see that half-graphs are connected and twin-free for all $k \geq 2$ and hence, are both FD- and \ftdadmis ~in this case.

%%%%%%%%%%%%%%%
%%%%%%%%%%%%%%%%

In \cite{FGRS_2021} it was shown that the only graphs whose OTD-numbers equal the order of the graph are the disjoint unions of half-graphs. In particular, we have $\otd(B_k) = 2k$. Combining this result with Corollary~\ref{cor_lb_PL-numbers} showing that $\otd(G) \leq \ftd(G)$ holds for all graphs and Theorem \ref{thm_pd-ptd_relations} showing $\fd(G + K_1) = \ftd(G) + 1$ together yields:
%%%%%%%%%%%%%%%%%%%%%%%%%
\begin{corollary}\label{cor_PL_half-graphs}
For a graph $G=(V,E)$ of order $n$ being the disjoint union of half-graphs, we have $\ftd(G) = n$ and $\fd(G + K_1) = \ftd(G) + 1 = n+1$. 
\end{corollary}
%%%%%%%%%%%%%%%%%%%%%%%%
Hence, there are graphs where the FD- and FTD-numbers equal the order of the graph.
A further example is $B_2=P_4$ with $\fd(P_4) = 4$.

\begin{theorem} \label{thm:half_graphs}
For a half-graph $B_k$ with $k \geq 3$, we have $\fd(B_k) = 2k-1$.
\end{theorem}

\begin{proof}
Consider a half-graph $B_k=(U \cup W, E)$ with vertices in $U = \{u_1, \ldots , u_k\}$, $W = \{w_1, \ldots , w_k \}$ and edges $u_iw_j$ if and only if $i \leq j$ and $k \geq 3$. 

In order to calculate $\fd(B_k)$, we first construct $\hyp_\FD(B_k)$. 
From the definition of FD-codes and Theorem~\ref{thm_PL-char}, we know that $\hyp_\FD(B_k)$ is composed of the closed neighborhoods
\begin{itemize}
\item $N[u_i] = \{u_i\} \cup \{w_i, \ldots , w_k\}$ for all $u_i \in U$,
\item $N[w_j] = \{u_1, \ldots , u_j\} \cup \{w_j\}$ for all $w_j \in W$, 
\end{itemize}
the symmetric differences of closed neighborhoods of adjacent vertices
\begin{itemize}
\item $N[u_i] \bigtriangleup N[w_j] = \{u_1, \ldots , u_{i-1}, u_{i+1} \ldots , u_j\} \cup \{w_i, \ldots , w_{j-1}, w_{j+1} \ldots , w_k\}$ for $i < j$,
\item $N[u_i] \bigtriangleup N[w_i] = \{u_1, \ldots , u_{i-1}\} \cup \{w_{i+1} \ldots , w_k\}$
\end{itemize}
and the symmetric differences of open neighborhoods of non-adjacent vertices
\begin{itemize}
\item $N(u_i) \bigtriangleup N(u_j) = \{w_i, \ldots , w_{j-1}\}$ for $i < j$,
\item $N(w_i) \bigtriangleup N(w_j) = \{u_{i+1} \ldots , u_j\}$ for $i < j$,
\item $N(u_i) \bigtriangleup N(w_j) = \{w_i, \ldots , w_k\} \cup \{u_1, \ldots , u_j\}$ for $i > j$.
\end{itemize}
In particular, we see that 
\begin{itemize}
\item $N(u_i) \bigtriangleup N(u_{i+1}) = \{w_i\}$ for $1 \leq i < k$,
\item $N(w_i) \bigtriangleup N(w_{i+1}) = \{u_{i+1}\}$ for $1 \leq i < k$, and 
\item $N(u_k) \bigtriangleup N(w_1) = \{w_k\} \cup \{u_1\}$
\end{itemize}
holds. Therefore, for $k \geq 3$, all neighborhoods and all other symmetric differences are redundant hyperedges (note that this is not the case for $k = 3$ as then $N[u_1] \bigtriangleup N[w_1] = \{w_2\}$ and $N[u_2] \bigtriangleup N[w_2] = \{u_1\}$ holds and makes $N(u_2)\bigtriangleup N(w_1) = \{w_2\} \cup \{u_1\}$ redundant).  
This clearly shows that $U \cup W \setminus \{w_k\}$ and $U \setminus \{u_1\} \cup W$ are the only two minimum FD-codes of the half-graph $B_k$ with $k \geq 3$, hence $\fd(B_k) = 2k-1$ indeed follows. 
\end{proof}

Hence, By Theorem~\ref{thm:half_graphs} and Corollary~\ref{cor_PL_half-graphs}, for half-graphs $B_k$ with $k \geq 3$, the FD- and FTD-numbers differ. This result combined with Theorem \ref{thm_pd-ptd_relations} shows that $G = B_k + K_1$ for $k \geq 3$ are examples of graphs where $\fd(G)$ is larger than the sum of the FD-numbers of its components. So far, this behavior was only observed for OD-numbers, see \cite{CW_ISCO2024}. Moreover, half-graphs are extremal graphs for the lower bounds $\otd(G) \leq \ftd(G)$ and $\max(\od(G), \otd(G)-1) \leq \fd(G)$ (the latter can be deduced by combining Corollary~\ref{cor_lb_PL-numbers} with the results from \cite{FGRS_2021} showing $\otd(B_k) = 2k$ and from \cite{CW_ISCO2024} showing $\od(B_k) = 2k-1$, respectively, for all $k \geq 1$).

\subsection{Headless spiders}

A graph $G=(U \cup W, E)$ is \emph{split} if its vertex set can be partitioned into a clique $Q$ and a stable set $S$. 
We next examine FD- and FTD-codes in two families of twin-free split graphs for which the exact X-numbers for $\X \in \{\LD, \LTD, \ID, \ITD, \OD, \OTD\}$ are known from \cite{ABLW_2022,ABW_2016,CW_ISCO2024,FL_2023}.

A \emph{headless spider} is a split graph with $Q = \{q_1, \ldots, q_k\}$ and $S = \{s_1, \ldots, s_k\}$. In addition, a headless spider is \emph{thin} (respectively, \emph{thick}) if $s_i$ is adjacent to $q_j$ if and only if $i = j$ (respectively, $i \neq j$). By definition, it is clear that the complement of a thin headless spider $H_k$ is a thick headless spider $\overline H_k$, and vice-versa. We have, for example, $H_2 = \overline H_2 = P_4$ and the thin headless spider $H_4$ is depicted in Figure~\ref{fig_exp_PL-codes}. Moreover, it is easy to check that thin and thick headless spiders have no twins.

It is known from \cite{ABLW_2022,CW_ISCO2024} that $\x(H_k) = k$ for $\X \in \{\LD, \LTD, \OD, \OTD\}$ and all $k \geq 3$. Furthermore, $\id(H_k) = k+1$ for $k \geq 3$ was shown in \cite{ABW_2016} and $\itd(H_k) = 2k-1$ for $k \geq 3$ in \cite{FL_2023}. 
We next analyze the X-numbers of thin headless spiders for X $\in$  \{FD, FTD\}.

\begin{theorem}\label{thm_thin-spider}
For a thin headless spider $H_k=(Q \cup S, E)$ with $k \geq 4$, we have $\fd(H_k) = 2k-2$ and $\ftd(H_k) = 2k-1$.
\end{theorem}

\begin{proof}
Consider a thin headless spider $H_k=(Q \cup S, E)$ with $Q = \{q_1, \ldots, q_k\}$, $S = \{s_1, \ldots, s_k\}$ and $k \geq 4$.
In order to calculate $\x(H_k)$ for X $\in$ \{FD, FTD\}, we first construct the two X-hypergraphs. From the definitions, Table \ref{tab_hypergraphs} and Theorem~\ref{thm_PL-char}, we know that the following hyperedges are involved: the closed or open neighborhoods, that is, 
\begin{itemize}
\item $N[q_i] = Q \cup \{s_i\}$ or $N(q_i) = Q \setminus \{q_i\} \cup \{s_i\}$ for all $q_i \in Q$,
\item $N[s_i] = \{q_i,s_i\}$ or $N(s_i) = \{q_i\}$ for all $s_i \in S$, 
\end{itemize}
the symmetric differences of closed neighborhoods of adjacent vertices, that is, 
\begin{itemize}
\item $N[s_i] \bigtriangleup N[q_i] = Q \setminus \{q_i\}$ for $1 \leq i \leq k$,
\item $N[q_i] \bigtriangleup N[q_j] = \{s_i,s_j\}$ for $1 \leq i < j \leq k$,
\end{itemize}
and the symmetric differences of open neighborhoods of non-adjacent vertices, that is, 
\begin{itemize}
\item $N(s_i) \bigtriangleup N(q_j) = Q \setminus \{q_i,q_j\} \cup \{s_j\}$ for $i \neq j$,
\item $N(s_i) \bigtriangleup N(s_j) = \{q_i,q_j\}$ for $1 \leq i < j \leq k$.
\end{itemize}
$\hyp_\FD(H_k)$ is composed of the closed neighborhoods and all the symmetric differences. In particular, we see that 
\begin{itemize}
\item $N[s_i] = \{q_i,s_i\}$ for $1 \leq i \leq k$, 
\item $N[q_i] \bigtriangleup N[q_j] = \{s_i,s_j\}$ for $1 \leq i < j \leq k$,
\item $N(s_i) \bigtriangleup N(s_j) = \{q_i,q_j\}$ for $1 \leq i < j \leq k$
\end{itemize}
belong to $\hyp_\FD(H_k)$. 
Therefore, $N[q_i]$ for all $q_i \in Q$, $N[s_i] \bigtriangleup N[q_i]$ for $1 \leq i \leq k$, and $N(s_i) \bigtriangleup N(q_j)$ for $i \neq j$ are redundant for $k \geq 4$ (note that $N(s_i) \bigtriangleup N(q_j)$ is not redundant if $k=3$). 
This shows that $Q \setminus \{q_i\} \cup S \setminus \{s_j\}$ for $i \neq j$ are the minimum FD-codes of $H_k$, hence $\fd(H_k) = 2k-2$ follows.

$\hyp_\FTD(H_k)$ is composed of the open neighborhoods and all the symmetric differences. In particular, we have that 
\begin{itemize}
\item $N(s_i) = \{q_i\}$ for $1 \leq i \leq k$, 
\item $N[q_i] \bigtriangleup N[q_j] = \{s_i,s_j\}$ for $1 \leq i < j \leq k$
\end{itemize}
belong to $\hyp_\FTD(H_k)$. 
Therefore, $N(q_i)$ for all $q_i \in Q$, $N[s_i] \bigtriangleup N[q_i]$ for $1 \leq i \leq k$, $N(s_i) \bigtriangleup N(q_j)$ for $i \neq j$ and $N(s_i) \bigtriangleup N(s_j)$ for $1 \leq i < j \leq k$ are redundant (even for $k \geq 3$). We conclude that $Q \cup S \setminus \{s_i\}$ for $1 \leq i \leq k$ are the minimum FTD-codes of $H_k$ and $\ftd(H_k) = 2k-1$ holds accordingly.
\end{proof}

It is known from \cite{ABLW_2022,CW_ISCO2024} that $\ld(\overline H_k) = \ltd(\overline H_k) = k-1$ and $\od(\overline H_k) = \otd(\overline H_k) = k+1$ holds for all $k \geq 3$. 
Furthermore, $\id(\overline H_k) = k$ for $k \geq 3$ was shown in \cite{ABW_2016}. 
We next determine the X-numbers of thick headless spiders for $\X \in  \{\FD, \FTD\}$. Moreover, for $\X \in \codes$, since the \ITD-number of thick headless spiders would be the only other \X-number that would be left to establish, we also integrate the said result in our next theorem.

\begin{theorem}\label{thm_thick-spider}
For a thick headless spider $\overline H_k=(Q \cup S, E)$ with $k \geq 4$, we have $\fd(\overline H_k) = \ftd(\overline H_k) = 2k-2$ and $\itd(\overline H_k) = k+1$.
\end{theorem}

\begin{proof}
Consider a thick headless spider $\overline H_k=(Q \cup S, E)$ with $Q = \{q_1, \ldots, q_k\}$, $S = \{s_1, \ldots, s_k\}$ and $k \geq 4$.
In order to calculate $\x(\overline H_k)$ for X $\in$ \{\FD, \FTD, \ITD\}, we first construct the three X-hypergraphs. From the definitions, Table \ref{tab_hypergraphs} and Theorem~\ref{thm_PL-char}, we know that the following hyperedges are involved: 
the closed or open neighborhoods, that is, 
\begin{itemize}
\item $N[s_i] = Q \setminus \{q_i\} \cup \{s_i\}$ or $N(s_i) = Q \setminus \{q_i\}$ for all $s_i \in S$,
\item $N[q_i] = Q \cup S \setminus \{s_i\}$ or $N(q_i) = Q \setminus \{q_i\} \cup S \setminus \{s_i\}$ for all $q_i \in Q$, 
\end{itemize}
the symmetric differences of closed neighborhoods of adjacent vertices, that is, 
\begin{itemize}
\item $N[s_i] \bigtriangleup N[q_j] = \{q_i\} \cup S \setminus \{s_i,s_j\}$ for $i \neq j$,
\item $N[q_i] \bigtriangleup N[q_j] = \{s_i,s_j\}$ for $i \neq j$,
\end{itemize}
and the symmetric differences of closed or open neighborhoods of non-adjacent vertices, that is, 
\begin{itemize}
\item $N[s_i] \bigtriangleup N[q_i] = \{q_i\} \cup S$ or $N(s_i) \bigtriangleup N(q_i) = S \setminus \{s_i\}$ for $1 \leq i \leq k$, 
\item $N[s_i] \bigtriangleup N[s_j] = \{q_i,q_j\} \cup \{s_i,s_j\}$ or $N(s_i) \bigtriangleup N(s_j) = \{q_i,q_j\}$ for $i \neq j$.
\end{itemize}
$\hyp_\FD(\overline H_k)$ is composed of the closed neighborhoods and the symmetric differences of closed neighborhoods of adjacent vertices as well as of open neighborhoods of non-adjacent vertices. In particular, we see that 
\begin{itemize}
\item $N[q_i] \bigtriangleup N[q_j] = \{s_i,s_j\}$ for $1 \leq i < j \leq k$,
\item $N(s_i) \bigtriangleup N(s_j) = \{q_i,q_j\}$ for $1 \leq i < j \leq k$
\end{itemize}
belong to $\hyp_\FD(\overline H_k)$. 
Therefore, all neighborhoods and the symmetric differences 
$N[s_i] \bigtriangleup N[q_j]$ for $i \neq j$ and $N(s_i) \bigtriangleup N(q_i)$ for $1 \leq i \leq k$ are redundant for $k \geq 4$ (note that $N[s_i] \bigtriangleup N[q_j]$ are not redundant if $k=3$). 
This shows that $Q \setminus \{q_i\} \cup S \setminus \{s_j\}$ for $1 \leq i,j \leq k$ are the minimum FD-codes of $\overline H_k$, hence $\fd(\overline H_k) = 2k-2$ follows.

$\hyp_\FTD(\overline H_k)$ is composed of the open neighborhoods and the symmetric differences of closed neighborhoods of adjacent vertices as well as of open neighborhoods of non-adjacent vertices. Again, we have that 
\begin{itemize}
\item $N[q_i] \bigtriangleup N[q_j] = \{s_i,s_j\}$ for $1 \leq i < j \leq k$,
\item $N(s_i) \bigtriangleup N(s_j) = \{q_i,q_j\}$ for $1 \leq i < j \leq k$
\end{itemize}
belong to $\hyp_\FTD(\overline H_k)$. 
Therefore, all neighborhoods and the other symmetric differences are redundant for $k \geq 4$.  
We conclude that $Q \setminus \{q_i\} \cup S \setminus \{s_j\}$ for $1 \leq i,j \leq k$ are also the minimum FTD-codes of $\overline H_k$ and thus also $\ftd(\overline H_k) = 2k-2$ holds.

$\hyp_\ITD(\overline H_k)$ is composed of the open neighborhoods and the symmetric differences of closed neighborhoods of distinct vertices. 
In particular, we see that 
\begin{itemize}
\item $N(s_i) = Q \setminus \{q_i\}$ for $1 \leq i \leq k$, 
\item $N[q_i] \bigtriangleup N[q_j] = \{s_i,s_j\}$ for $1 \leq i < j \leq k$
\end{itemize}
belong to $\hyp_\ITD(\overline H_k)$. 
Thus, $N(q_i)$ for all $q_i \in Q$ and all other symmetric differences are redundant for $k \geq 4$ (note that $N[s_i] \bigtriangleup N[q_j]$ are not redundant if $k=3$). 
This implies that two vertices from $Q$ and all but one vertices from $S$ have to be taken in every minimum ITD-code of $\overline H_k$. Thus, we have $\itd(\overline H_k) = k+1$.
\end{proof}

Hence, thin headless spiders are extremal graphs for the lower bounds $\itd(G) \leq \ftd(G)$ and $\itd(G) -1 \leq \fd(G)$, whereas thick headless spiders illustrate that the gap between the FD- or FTD-numbers and their lower bounds by Corollary \ref{cor_all_relations} can be large.

%%%%%%%%%%%%%%%%%%%%%%%%%%%%%%%%%%%

\section{Concluding remarks}

In this paper, we studied problems in the domain of identification problems. The contribution of this paper is to introduce a new separation property, called full-separation, and to study it in combination with both domination and total-domination resulting in FD-codes and FTD-codes, respectively. It turns out that FTD-codes have been introduced in the literature of identification problems under a differently formulated definition and under the name of SID-codes. In this paper, we characterized full-separating sets in different ways to show that they combine the requirements of both closed- and open-separation which justifies in particular how the full-separation property models monitoring systems that tackle the two studied fault types simultaneously. We addressed questions concerning the existence of FD- and FTD-codes in graphs, the relations between all X-numbers and particularly between the FD- and FTD-numbers.
We further showed that deciding if the FD- or FTD-number of a graph is below a given value is NP-complete. Moreover, it is in general NP-hard to decide if the two numbers differ on a graph. 
 We also provided some extremal examples of graphs with respect to lower and upper bounds on the FD- and FTD-numbers.
 
Our lines of future research include to pursue these studies of the FD- or FTD-codes by involving more graph families, for example, to find more examples of graphs where the FD- or FTD-numbers are not close to the order of the graph. In general it would be interesting to identify or even characterize the extremal cases, that is, graphs for which the FD- or FTD-numbers equal the upper bound (the order of the graph) or one of its lower bounds in terms of other X-numbers.

Finally, it would be interesting to study a broader set of problem variants in the area of domination / separation from a hypergraph point of view towards a more general classification of such problems. In this context, it is clear that the hypergraphs encoding domination (respectively, total-domination) are only composed of closed (respectively, open) neighborhoods, whereas hypergraphs encoding only the separation properties only contain the corresponding symmetric differences of neighborhoods for all pairs of vertices. Here, it would be interesting to see if other restrictions (for example, only to $\ddelta_a[G]$ or only to $N(G)$ and $\ddelta_a(G)$) encode further studied problems from the literature and then similar relations as in Figure~\ref{fig_relations} can be deduced to a large set of graph parameters.

%%%%%%%%%%%%%%%%%%%%%%%%%%%%%%%%%%%

\section*{Acknowledgements}
\noindent This research was financed by a public grant overseen by the French National Research Agency as part of the ``Investissements d’Avenir'' through the IMobS3 Laboratory of Excellence (ANR-10-LABX-0016), by the French government IDEX-ISITE initiative 16-IDEX-0001 (CAP 20-25) and the International Research Center ``Innovation Transportation and Production Systems'' of the I-SITE CAP 20-25.

\bibliographystyle{splncs04}

\end{document}